%
%
\font\hd=cmbx10 scaled\magstep1



\ifx\begin\undefined\else\global\advance\srcdepth by
1\expandafter \fi

\def\begin{}
\newcount\srcdepth
\srcdepth=1

\outer\def\bye{\global\advance\srcdepth by -1
  \ifnum\srcdepth=0
    \def\endcmd{\vfill\supereject\nopagenumbers\par\vfill\supereject\end}
  \else\def\endcmd{}\fi
  \endcmd
}



\def\initialize#1#2#3#4#5#6{
  \ifnum\srcdepth=1
  \magnification=#1
  \hsize = #2
  \vsize = #3
  \hoffset=#4
  \advance\hoffset by -\hsize
  \divide\hoffset by 2
  \advance\hoffset by -1truein
  \voffset=#5
  \advance\voffset by -\vsize
  \divide\voffset by 2
  \advance\voffset by -1truein
  \advance\voffset by #6
  \baselineskip=13pt
  \emergencystretch = 0.05\hsize
  \parskip=3pt plus1pt minus.5pt
  \fi
}

\def\print{\initialize{1095}
  {6.5truein}{8.5truein}{8.5truein}{11truein}{-.0625truein}}

\overfullrule=0pt

\newif\ifblackboardbold

\blackboardboldtrue


\font\sectionfont=cmbx12

\font\scriptit=cmti10 at 7pt
\font\scriptsl=cmsl10 at 7pt
\scriptfont\itfam=\scriptit
\scriptfont\slfam=\scriptsl


\newfam\msamfam  
\font\tenmsam=msam10
\font\sevenmsam=msam7
\font\fivemsam=msam5
\textfont\msamfam=\tenmsam
\scriptfont\msamfam=\sevenmsam
\scriptscriptfont\msamfam=\fivemsam

\newfam\msbmfam
\font\tenmsbm=msam10
\font\sevenmsbm=msam7
\font\fivemsbm=msam5
\textfont\msbmfam=\tenmsbm
\scriptfont\msbmfam=\sevenmsbm
\scriptscriptfont\msbmfam=\fivemsbm

\newfam\eufmfam  
\font\teneufm=eufm10
\font\seveneufm=eufm7
\font\fiveeufm=eufm5
\textfont\eufmfam\teneufm
\scriptfont\eufmfam\seveneufm
\scriptscriptfont\eufmfam\fiveeufm

\newfam\bboldfam
\ifblackboardbold
\font\tenbbold=msbm10
\font\sevenbbold=msbm7
\font\fivebbold=msbm5
\textfont\bboldfam=\tenbbold
\scriptfont\bboldfam=\sevenbbold
\scriptscriptfont\bboldfam=\fivebbold
\def\bbold{\fam\bboldfam\tenbbold}
\else
\def\bbold{\bf}
\fi
\newcount\amsfamcount 
\newcount\classcount   
\newcount\positioncount
\newcount\codecount
\newcount\n             
\def\newsymbol#1#2#3#4#5{               
\n="#2                                  
\ifnum\n=1 \amsfamcount=\msamfam\else   
\ifnum\n=2 \amsfamcount=\msbmfam\else   
\ifnum\n=3 \amsfamcount=\eufmfam
\fi\fi\fi
\multiply\amsfamcount by "100           
\classcount="#3                 
\multiply\classcount by "1000           
\positioncount="#4#5            
\codecount=\classcount                  
\advance\codecount by \amsfamcount      
\advance\codecount by \positioncount
\mathchardef#1=\codecount}              


\ifblackboardbold

\fi








\newlinechar=`@
\def\forwardmsg#1#2#3{\immediate\write16{@*!*!*!* forward reference should
be: @\noexpand\forward{#1}{#2}{#3}@}}
\def\nodefmsg#1{\immediate\write16{@*!*!*!* #1 is an undefined reference@}}

\def\forwardsub#1#2{\def\newref{{#2}{#1}}}

\def\forward#1#2#3{%
\expandafter\expandafter\expandafter\forwardsub\expandafter{#3}{#2}
\expandafter\ifx\csname#1\endcsname\relax\else%
\expandafter\ifx\csname#1\endcsname\newref\else%
\forwardmsg{#1}{#2}{#3}\fi\fi%
\expandafter\let\csname#1\endcsname\newref}

\def\firstarg#1{\expandafter\argone #1}\def\argone#1#2{#1}
\def\secondarg#1{\expandafter\argtwo #1}\def\argtwo#1#2{#2}

\def\ref#1{\expandafter\ifx\csname#1\endcsname\relax
  {\nodefmsg{#1}\bf`#1'}\else
  \expandafter\firstarg\csname#1\endcsname
  ~\expandafter\secondarg\csname#1\endcsname\fi}

\def\refs#1{\expandafter\ifx\csname#1\endcsname\relax
  {\nodefmsg{#1}\bf`#1'}\else
  \expandafter\firstarg\csname #1\endcsname
  s~\expandafter\secondarg\csname#1\endcsname\fi}

\def\refn#1{\expandafter\ifx\csname#1\endcsname\relax
  {\nodefmsg{#1}\bf`#1'}\else
  \expandafter\secondarg\csname #1\endcsname\fi}



\def\widow#1{\vskip 0pt plus#1\vsize\goodbreak\vskip 0pt plus-#1\vsize}



\def\marginlabel#1{}

\def\showlabelsabove{
\font\labelfont=cmss10 at 6pt
\def\marginlabel##1{\rlap{\smash{\raise 10pt\hbox{\labelfont##1}}}}
}

\newcount\seccount
\newcount\proccount
\seccount=0
\proccount=0

\def\stdskip{\vskip 9pt plus3pt minus 3pt}
\def\stdbreak{\par\removelastskip\penalty-100\stdskip}

\def\proof{\stdbreak\noindent{\sl Proof }}

\def\qed{\vrule height 1.2ex width .9ex depth .1ex}

\def\Box{
  \ifmmode\eqno\qed
  \else\ifvmode\removelastskip\line{\hfil\qed}
  \else\unskip\quad\hskip-\hsize
    \hbox{}\hskip\hsize minus 1em\qed\par
  \fi\stdbreak\fi}

\def\references{
  \removelastskip
  \widow{.05}
  \vskip 24pt plus 6pt minus 6 pt
  \leftline{\sectionfont References}
  \nobreak\stdskip\noindent}

\def\ifempty#1#2\endB{\ifx#1\endA}
\def\makeref#1#2#3{\ifempty#1\endA\endB\else\forward{#1}{#2}{#3}\fi}

\outer\def\section#1 #2\par{
  \removelastskip
  \global\advance\seccount by 1
  \global\proccount=0\relax
                \edef\numtoks{\number\seccount}
  \makeref{#1}{Section}{\numtoks}
  \widow{.05}
  \vskip 24pt plus 6pt minus 6 pt
  \message{#2}
  \leftline{\marginlabel{#1}\sectionfont\numtoks\quad #2}
  \nobreak\stdskip}

\def\proclamation#1#2{
  \outer\expandafter\def\csname#1\endcsname##1 ##2\par{
  \stdbreak
  \global\advance\proccount by 1
  \edef\numtoks{\number\seccount.\number\proccount}
  \makeref{##1}{#2}{\numtoks}
  \noindent{\marginlabel{##1}\bf #2 \numtoks\enspace}
  {\sl##2\par}
  \stdbreak}}

\def\othernumbered#1#2{
  \outer\expandafter\def\csname#1\endcsname##1{
  \stdbreak
  \global\advance\proccount by 1
  \edef\numtoks{\number\seccount.\number\proccount}
  \makeref{##1}{#2}{\numtoks}
  \noindent{\marginlabel{##1}\bf #2 \numtoks\enspace}}}

\proclamation{definition}{Definition}
\proclamation{lemma}{Lemma}
\proclamation{proposition}{Proposition}
\proclamation{theorem}{Theorem}
\proclamation{corollary}{Corollary}
\proclamation{conjecture}{Conjecture}

\othernumbered{example}{Example}
\othernumbered{remark}{Remark}
\othernumbered{construction}{Construction}


\def\figure#1{
 \global\advance\figcount by 1
 \goodbreak
 \midinsert#1\smallskip
 \centerline{Figure~\number\figcount}
 \endinsert}

\def\capfigure#1#2{
 \global\advance\figcount by 1
 \goodbreak
 \midinsert#1\smallskip
 \vbox{\small\noindent {\bf Figure~\number\figcount:} #2}
 \endinsert}

\def\capfigurepair#1#2#3#4{
 \goodbreak
 \midinsert
 #1\smallskip
 \global\advance\figcount by 1
 \vbox{\small\noindent {\bf Figure~\number\figcount:} #2}
 \vskip 12pt
 #3\smallskip
 \global\advance\figcount by 1
 \vbox{\small\noindent {\bf Figure~\number\figcount:} #4}
 \endinsert}


\def\baretable#1#2{
\vbox{\offinterlineskip\halign{
 \strut\kern #1\hfil##\kern #1
 &&\kern #1\hfil##\kern #1\cr
 #2
}}}

\def\gridtablesub#1#2#3{
\vbox{\offinterlineskip\halign{
 \strut\vrule\kern #1\hfil##\hfil\kern #2\vrule
 &&\kern #1\hfil##\kern #2\vrule\cr
 \noalign{\hrule}
 #3
 \noalign{\hrule}
}}}





\input epsf

\newcount\figcount
\figcount=0
\newbox\drawing
\newcount\drawbp
\newdimen\drawx
\newdimen\drawy
\newdimen\ngap
\newdimen\sgap
\newdimen\wgap
\newdimen\egap

\def\drawbox#1#2#3{\vbox{
  \setbox\drawing=\vbox{\offinterlineskip\epsfbox{#2.eps}\kern 0pt}
  \drawbp=\epsfurx
  \advance\drawbp by-\epsfllx\relax
  \multiply\drawbp by #1
  \divide\drawbp by 100
  \drawx=\drawbp truebp
  \ifdim\drawx>\hsize\drawx=\hsize\fi
  \epsfxsize=\drawx
  \setbox\drawing=\vbox{\offinterlineskip\epsfbox{#2.eps}\kern 0pt}
  \drawx=\wd\drawing
  \drawy=\ht\drawing
  \ngap=0pt \sgap=0pt \wgap=0pt \egap=0pt
  \setbox0=\vbox{\offinterlineskip
    \box\drawing \ifgridlines\drawgrid\drawx\drawy\fi #3}
  \kern\ngap\hbox{\kern\wgap\box0\kern\egap}\kern\sgap}}

\def\draw#1#2#3{
  \setbox\drawing=\drawbox{#1}{#2}{#3}
  \global\advance\figcount by 1
  \goodbreak
  \midinsert
  \centerline{\ifgridlines\boxgrid\drawing\fi\box\drawing}
  \smallskip
  \vbox{\offinterlineskip
    \centerline{Figure~\number\figcount}
    \smash{\marginlabel{#2}}}
  \endinsert}

\def\capdraw#1#2#3#4{
  \setbox\drawing=\drawbox{#1}{#2}{#3}
  \global\advance\figcount by 1
  \goodbreak
  \midinsert
  \centerline{\ifgridlines\boxgrid\drawing\fi\box\drawing}
  \smallskip
  \vbox{\offinterlineskip
    \vskip 4pt
    \vbox{\lineskip=3pt\small\noindent {\bf Figure~\number\figcount:} #4}
    \smash{\marginlabel{#2}}}
  \endinsert}

\def\capdrawpair#1#2#3#4#5#6#7#8{
  \goodbreak
  \midinsert
  \setbox\drawing=\drawbox{#1}{#2}{#3}
  \global\advance\figcount by 1
  \centerline{\ifgridlines\boxgrid\drawing\fi\box\drawing}
  \smallskip
  \vbox{\offinterlineskip
    \vskip 4pt
    \vbox{\lineskip=3pt\small\noindent {\bf Figure~\number\figcount:} #4}
    \smash{\marginlabel{#2}}}
  \vskip 12pt
  \setbox\drawing=\drawbox{#5}{#6}{#7}
  \global\advance\figcount by 1
  \centerline{\ifgridlines\boxgrid\drawing\fi\box\drawing}
  \smallskip
  \vbox{\offinterlineskip
    \vskip 4pt
    \vbox{\lineskip=3pt\small\noindent {\bf Figure~\number\figcount:} #8}
    \smash{\marginlabel{#6}}}
  \endinsert}

\def\drawnoname#1#2#3{
  \setbox\drawing=\drawbox{#1}{#2}{#3}
  \global\advance\figcount by 1
  \goodbreak
  \midinsert
  \centerline{\ifgridlines\boxgrid\drawing\fi\box\drawing}
  \smallskip
  \endinsert}

\def\nextfigtoks{%
  \advance\figcount by 1%
  \edef\numtoks{\number\figcount}%
  \advance\figcount by -1}

\newif\ifgridlines
\newbox\figtbox
\newbox\figgbox
\newdimen\figtx
\newdimen\figty

\newdimen\bwd
\bwd=2sp 

\def\hline#1{\vbox{\smash{\hbox to #1{\leaders\hrule height \bwd\hfil}}}}

\def\vline#1{\hbox to 0pt{%
  \hss\vbox to #1{\leaders\vrule width \bwd\vfil}\hss}}

\def\clap#1{\hbox to 0pt{\hss#1\hss}}
\def\vclap#1{\vbox to 0pt{\offinterlineskip\vss#1\vss}}

\def\hstutter#1#2{\hbox{%
  \setbox0=\hbox{#1}%
  \hbox to #2\wd0{\leaders\box0\hfil}}}

\def\vstutter#1#2{\vbox{
  \setbox0=\vbox{\offinterlineskip #1}
  \dp0=0pt
  \vbox to #2\ht0{\leaders\box0\vfil}}}

\def\crosshairs#1#2{
  \dimen1=.002\drawx
  \dimen2=.002\drawy
  \ifdim\dimen1<\dimen2\dimen3\dimen1\else\dimen3\dimen2\fi
  \setbox1=\vclap{\vline{2\dimen3}}
  \setbox2=\clap{\hline{2\dimen3}}
  \setbox3=\hstutter{\kern\dimen1\box1}{4}
  \setbox4=\vstutter{\kern\dimen2\box2}{4}
  \setbox1=\vclap{\vline{4\dimen3}}
  \setbox2=\clap{\hline{4\dimen3}}
  \setbox5=\clap{\copy1\hstutter{\box3\kern\dimen1\box1}{6}}
  \setbox6=\vclap{\copy2\vstutter{\box4\kern\dimen2\box2}{6}}
  \setbox1=\vbox{\offinterlineskip\box5\box6}
  \smash{\vbox to #2{\hbox to #1{\hss\box1}\vss}}}

\def\boxgrid#1{\rlap{\vbox{\offinterlineskip
  \setbox0=\hline{\wd#1}
  \setbox1=\vline{\ht#1}
  \smash{\vbox to \ht#1{\offinterlineskip\copy0\vfil\box0}}
  \smash{\vbox{\hbox to \wd#1{\copy1\hfil\box1}}}}}}

\def\drawgrid#1#2{\vbox{\offinterlineskip
  \dimen0=\drawx
  \dimen1=\drawy
  \divide\dimen0 by 10
  \divide\dimen1 by 10
  \setbox0=\hline\drawx
  \setbox1=\vline\drawy
  \smash{\vbox{\offinterlineskip
    \copy0\vstutter{\kern\dimen1\box0}{10}}}
  \smash{\hbox{\copy1\hstutter{\kern\dimen0\box1}{10}}}}}

\def\figtext#1#2#3#4#5{
  \setbox\figtbox=\vbox{\hbox{#5}\kern 0pt}
  \figtx=-#3\wd\figtbox \figty=-#4\ht\figtbox
  \advance\figtx by #1\drawx \advance\figty by #2\drawy
  \dimen0=\figtx \advance\dimen0 by\wd\figtbox \advance\dimen0 by-\drawx
  \ifdim\dimen0>\egap\global\egap=\dimen0\fi
  \dimen0=\figty \advance\dimen0 by\ht\figtbox \advance\dimen0 by-\drawy
  \ifdim\dimen0>\ngap\global\ngap=\dimen0\fi
  \dimen0=-\figtx
  \ifdim\dimen0>\wgap\global\wgap=\dimen0\fi
  \dimen0=-\figty
  \ifdim\dimen0>\sgap\global\sgap=\dimen0\fi
  \smash{\rlap{\vbox{\offinterlineskip
    \hbox{\hbox to \figtx{}\ifgridlines\boxgrid\figtbox\fi\box\figtbox}
    \vbox to \figty{}
    \ifgridlines\crosshairs{#1\drawx}{#2\drawy}\fi
    \kern 0pt}}}}


\def\hpad#1#2#3{\hbox{\kern #1\hbox{#3}\kern #2}}
\def\vpad#1#2#3{\setbox0=\hbox{#3}\vbox{\kern #1\box0\kern #2}}




\def\stack#1#2#3{\vbox{\offinterlineskip
  \setbox2=\hbox{#2}
  \setbox3=\hbox{#3}
  \dimen0=\ifdim\wd2>\wd3\wd2\else\wd3\fi
  \hbox to \dimen0{\hss\box2\hss}
  \kern #1
  \hbox to \dimen0{\hss\box3\hss}}}


\def\hexp#1{%
  \setbox0=\hbox{${}^{#1}$}%
  \hbox to .5\wd0{\box0\hss}}

\def\hsub#1{%
  \setbox0=\hbox{${}_{#1}$}%
  \hbox to .5\wd0{\box0\hss}}



\def\bmatrix#1#2{{\left[\vcenter{\halign
  {&\kern#1\hfil$##\mathstrut$\kern#1\cr#2}}\right]}}

\def\rightarrowmat#1#2#3{
  \setbox1=\hbox{\small\kern#2$\bmatrix{#1}{#3}$\kern#2}
  \,\vbox{\offinterlineskip\hbox to\wd1{\hfil\copy1\hfil}
    \kern 3pt\hbox to\wd1{\rightarrowfill}}\,}

\def\leftarrowmat#1#2#3{
  \setbox1=\hbox{\small\kern#2$\bmatrix{#1}{#3}$\kern#2}
  \,\vbox{\offinterlineskip\hbox to\wd1{\hfil\copy1\hfil}
    \kern 3pt\hbox to\wd1{\leftarrowfill}}\,}

\def\rightarrowbox#1#2{
  \setbox1=\hbox{\kern#1\hbox{\small #2}\kern#1}
  \,\vbox{\offinterlineskip\hbox to\wd1{\hfil\copy1\hfil}
    \kern 3pt\hbox to\wd1{\rightarrowfill}}\,}

\def\leftarrowbox#1#2{
  \setbox1=\hbox{\kern#1\hbox{\small #2}\kern#1}
  \,\vbox{\offinterlineskip\hbox to\wd1{\hfil\copy1\hfil}
    \kern 3pt\hbox to\wd1{\leftarrowfill}}\,}








\def\quiremacro#1#2#3#4#5#6#7#8#9{
  \expandafter\def\csname#1\endcsname##1{
  \ifnum\srcdepth=1
  \magnification=#2
  \input quire
  \hsize=#3
  \vsize=#4
  \htotal=#5
  \vtotal=#6
  \shstaplewidth=#7
  \shstaplelength=#8
  \hoffset=\htotal
  \advance\hoffset by -\hsize
  \divide\hoffset by 2
  \ifnum\vsize<\vtotal
    \voffset=\vtotal
    \advance\voffset by -\vsize
    \divide\voffset by 2
  \fi
  \advance\voffset by #9
  \shhtotal=2\htotal
  \baselineskip=13pt
  \emergencystretch = 0.05\hsize
  \horigin=0.0truein
  \vorigin=0.0truein
  \shthickness=0pt
  \shoutline=0pt
  \shcrop=0pt
  \shvoffset=-1.0truein
  \ifnum##1>0\quire{#1}\else\qtwopages\fi
  \fi
}}



\quiremacro{letterbooklet} 
{1000}{4.79452truein}{7truein}{5.5truein}{8.5truein}{0.01pt}{0.66truein}{-.0625t
ruein}

\quiremacro{Afourbooklet}
{1095}{5.25truein}{7truein}{421truept}{595truept}{0.01pt}{0.66truein}{-.0625true
in}

\quiremacro{legalbooklet}
{1095}{5.25truein}{7truein}{7.0truein}{8.5truein}{0.01pt}{0.66truein}{-.0625true
in}

\quiremacro{twoupsub} 
{895}{4.5truein}{7truein}{5.5truein}{8.5truein}{0pt}{0pt}{.0625truein}


\quiremacro{Afourviewsub} 
{1000}{5.0228311in}{7.7625571in}{421truept}{595truept}{0.1pt}{0.5\vtotal}{-.0625
truein}


\quiremacro{viewsub}
{1095}{5.5truein}{8.5truein}{461truept}{666truept}{0.1pt}{0.5\vtotal}{-.125truei
n}


\newcount\countA
\newcount\countB
\newcount\countC

\def\monthname{\begingroup
  \ifcase\number\month
    \or January\or February\or March\or April\or May\or June\or
    July\or August\or September\or October\or November\or December\fi
\endgroup}

\def\dayname{\begingroup
  \countA=\number\day
  \countB=\number\year
  \advance\countA by 0 
  \advance\countA by \ifcase\month\or
    0\or 31\or 59\or 90\or 120\or 151\or
    181\or 212\or 243\or 273\or 304\or 334\fi
  \advance\countB by -1995
  \multiply\countB by 365
  \advance\countA by \countB
  \countB=\countA
  \divide\countB by 7
  \multiply\countB by 7
  \advance\countA by -\countB
  \advance\countA by 1
  \ifcase\countA\or Sunday\or Monday\or Tuesday\or Wednesday\or
    Thursday\or Friday\or Saturday\fi
\endgroup}

\def\timename{\begingroup
   \countA = \time
   \divide\countA by 60
   \countB = \countA
   \countC = \time
   \multiply\countA by 60
   \advance\countC by -\countA
   \ifnum\countC<10\toks1={0}\else\toks1={}\fi
   \ifnum\countB<12 \toks0={\sevenrm AM}
     \else\toks0={\sevenrm PM}\advance\countB by -12\fi
   \relax\ifnum\countB=0\countB=12\fi
   \hbox{\the\countB:\the\toks1 \the\countC \thinspace \the\toks0}
\endgroup}

\def\timestamp{\dayname, \the\day\ \monthname\ \the\year, \timename}


\print


\def\enma#1{{\ifmmode#1\else$#1$\fi}}

\def\mathbb#1{{\bbold #1}}
\def\mathbf#1{{\bf #1}}


\def\PP{\enma{\mathbb{P}}}



\def\CC{\enma{\mathbf{C}}}


\def\setdef#1#2{\enma{\{\;#1\;\,|\allowbreak
  \;\,#2\;\}}}


\def\im{\mathop{\rm im}\nolimits}

\newsymbol\boxtimes1202

\font\tenmsam=msam10
\font\sevenmsam=msam7
\font\fivemsam=msam5
\newfam\msamfam         
\textfont\msamfam\tenmsam
\scriptfont\msamfam\sevenmsam
\scriptscriptfont\msamfam\fivemsam

%

\font\tenmsbm=msbm10
\font\sevenmsbm=msbm7
\font\fivemsbm=msbm5
\newfam\msbmfam         
\textfont\msbmfam\tenmsbm
\scriptfont\msbmfam\sevenmsbm
\scriptscriptfont\msbmfam\fivemsbm

%

\font\teneufm=eufm10
\font\seveneufm=eufm7
\font\fiveeufm=eufm5
\newfam\eufmfam        
\textfont\eufmfam\teneufm
\scriptfont\eufmfam\seveneufm
\scriptscriptfont\eufmfam\fiveeufm

%

%
%
%
%
%
%
%
\newcount\amsfamcount 
\newcount\classcount   
\newcount\positioncount
\newcount\codecount
\newcount\n             
\def\newsymbol#1#2#3#4#5{               
\n="#2                                  
\ifnum\n=1 \amsfamcount=\msamfam\else   
\ifnum\n=2 \amsfamcount=\msbmfam\else   
\ifnum\n=3 \amsfamcount=\eufmfam
\fi\fi\fi
\multiply\amsfamcount by "100           
\classcount="#3                 
\multiply\classcount by "1000           
\positioncount="#4#5            
\codecount=\classcount                  
\advance\codecount by \amsfamcount      
\advance\codecount by \positioncount
\mathchardef#1=\codecount}              
\newcount\famcnt 
\newcount\classcnt   
\newcount\positioncnt
\newcount\codecnt
\def\newmathsymbol#1#2#3#4#5{          
\famcnt=#2                      
\multiply\famcnt by "100        
\classcnt="#3                   
\multiply\classcnt by "1000     
\positioncnt="#4#5              
\codecnt=\classcnt              
\advance\codecnt by \famcnt     
\advance\codecnt by \positioncnt
\mathchardef#1=\codecnt}        

%
\newsymbol\varnothing203F

\newsymbol\SEMI226F
\def\rtimes{\mathop{\SEMI}}
\newsymbol\smallsetminus2272
\def\setminus{\mathop{\smallsetminus}}

\def\Box{
  \ifmmode\eqno\qed
  \else\ifvmode\removelastskip\line{\hfil\qed}
  \else\unskip\quad\hskip-\hsize
    \hbox{}\hskip\hsize minus 1em\qed\par
  \fi\stdbreak\fi}

\def\PGL {\mathop{\rm PGL}\nolimits} 
\def\GL {\mathop{\rm GL}\nolimits} 

\def\PSL {\mathop{\rm PSL}\nolimits}
\def\Pic {\mathop{\rm Pic}\nolimits}
\def\Gr  {\mathop{\rm Gr}\nolimits}

\def\Proj {\mathop{\rm Proj}\nolimits}

\def\mod {\mathop{\rm mod}}

\def\diag {\mathop{\rm diag}\nolimits}
\def\mult {\mathop{\rm mult}\nolimits}
\def\vt{\enma{\vartheta}}
\def\boldz{{\bf Z}}
\def\HHH{{\bf H}}
\def \diag {\mathop{\rm diag}\nolimits}
\def\Pone{{\bf P}^1}
\def\Ptwo{{\bf P}^2}
\def\Pthree{{\bf P}^3}
\def\Pfour{{\bf P}^4}
\def\Pfive{{\bf P}^5}
\def\Psix{{\bf P}^6}

\def\A{{\cal A}}

\def\E{{\cal E}}

\def\G{{\cal G}}
\def\H{{\cal H}}
\def\I{{\cal I}}

\def\L{{\cal L}}
\def\O{{\cal O}}
\def\P{{\bf P}}

\def\X{{\cal X}}

\def\dual#1{{#1}^{\scriptscriptstyle \vee}}
\def\exact#1#2#3{0\longrightarrow#1\longrightarrow#2
\longrightarrow#3\longrightarrow0}

\def\tenpoint{%
\textfont0=\tenrm \scriptfont0=\sevenrm
\scriptscriptfont0=\fiverm \def\rm{\fam0\tenrm}%
\textfont1=\teni \scriptfont1=\seveni
\scriptscriptfont1=\fivei \def\oldstyle{\fam1\teni}%
\textfont2=\tensy \scriptfont2=\sevensy
\scriptscriptfont2=\fivesy
\textfont\itfam=\tenit \def\it{\fam\itfam\tenit}%
\textfont\slfam=\tensl \def\sl{\fam\slfam\tensl}%
\textfont\ttfam=\tentt \def\tt{\fam\tfam\tentt}%
\textfont\bffam=\tenbf \scriptfont\bffam=\sevenbf
\scriptscriptfont\bffam=\fivebf \def\bf{\fam\bffam\tenbf}%
\abovedisplayskip=12pt plus 3pt minus 9pt
\belowdisplayskip=\abovedisplayskip
\abovedisplayshortskip=0pt plus 3pt
\belowdisplayshortskip=7pt plus 3pt minus 4pt
\smallskipamount=3pt plus 1pt minus 1pt
\medskipamount=6pt plus 2pt minus 2pt
\bigskipamount=12pt plus 4pt minus 4pt
\setbox\strutbox=\hbox{\vrule height8.5pt depth3.5pt width0pt}%
\normalbaselineskip=12pt \normalbaselines \rm}

\def\eightpoint{%
\font\eightrm=cmr8%
\font\sixrm=cmr6%
\font\eighti=cmmi6%
\font\eightit=cmti8%
\font\sixi=cmmi6%
\font\sixit=cmti6%
\font\eightsy=cmsy8%
\font\sixsy=cmsy6%
\font\eightsl=cmsl8%
\font\eighttt=cmtt6%
\font\eightbf=cmbx8%
\font\sixbf=cmbx6%
\textfont0=\eightrm \scriptfont0=\sixrm
\scriptscriptfont0=\fiverm \def\rm{\fam0\eightrm}%
\textfont1=\eighti \scriptfont1=\sixi
\scriptscriptfont1=\fivei \def\oldstyle{\fam1\eighti}%
\textfont2=\eightsy \scriptfont2=\sixsy
\scriptscriptfont2=\fivesy
\textfont\itfam=\eightit \def\it{\fam\itfam\eightit}%
\textfont\slfam=\eightsl \def\sl{\fam\slfam\eightsl}%
\textfont\ttfam=\eighttt \def\tt{\fam\tfam\eighttt}%
\textfont\bffam=\eightbf \scriptfont\bffam=\sixbf
\scriptscriptfont\bffam=\fivebf \def\bf{\fam\bffam\eightbf}%
\abovedisplayskip=9pt plus 3pt minus 9pt
\belowdisplayskip=\abovedisplayskip
\abovedisplayshortskip=0pt plus 2pt
\belowdisplayshortskip=5pt plus 2pt minus 3pt
\smallskipamount=2pt plus 1pt minus 1pt
\medskipamount=4pt plus 2pt minus 2pt
\bigskipamount=9pt plus 4pt minus 4pt
\setbox\strutbox=\hbox{\vrule height 7pt depth 2pt width 0pt}%
\normalbaselineskip=9pt \normalbaselines \rm}

\def\rto{\raise.5ex\hbox{$\scriptscriptstyle ---\!\!\!>$}}


\forward{1.7}{Section}{5}
\forward{defect6}{}{}
\forward {defect8}{Remark}{6.9}
\forward {rat.map6}{Remark}{6.9}

\centerline{\hd Calabi-Yau Three-folds and Moduli of Abelian Surfaces II}
\medskip
\centerline{\it Mark Gross\footnote{*}{Supported by NSF grant DMS-0805328}
and Sorin Popescu\footnote{**}{Partially supported by
NSF grant DMS-0502070, DMS-0083361, and MSRI, Berkeley. Current Address:
Renaissance Technologies, 600 Route 25A, East Setauket, NY 11733}}
\medskip

\centerline{1 August, 2009}
\medskip
{\settabs 3 \columns
\+UCSD Mathematics&&
Department of Mathematics\cr
\+9500 Gilman Drive&&
Stony Brook University\cr
\+La Jolla, CA 92093-0112&&
Stony Brook, NY 11794-3651\cr
\+mgross@math.ucsd.edu&&
sorin@rentec.com\cr}

\bigskip
\bigskip

{\hd \S 0. Introduction.}

The main goal of this paper, which is a continuation of [GP1], [GP2] and [GP3],
is to describe birational models for moduli spaces $\A_{d}$
of polarized abelian surfaces of type $(1,d)$ for small values of $d$,
and for moduli spaces of such polarized abelian surfaces with 
suitably defined partial or canonical level structure. 
We can then decide the uniruledness, unirationality or rationality
of nonsingular models of compactifications of these moduli spaces
(which are quasi-projective $3$-folds, possibly singular).
Since these properties are birational invariants, this determines the
corresponding properties of those moduli spaces.

We will use in the sequel definitions and notation as in [GP1], [GP2], [GP3]; 
see also [Mu1], [LB] and [HKW] for basic facts concerning 
abelian varieties and their moduli, as well as
for the definition of a canonical level structure. We will also
make use of  partial level structures which are introduced and
described in \S 1.

Using a version of the Maass-Kurokawa lifting, Gritsenko [Gri1], [Gri2]
proved the existence of weight 3 cusp forms with
respect to the paramodular group $\Gamma_d$, for almost
all values of $d$. Since one knows the dimension for
the space of Jacobi cusp forms one deduces lower
bounds for the dimensions of the spaces of cusp 
forms with respect to the paramodular group $\Gamma_d$.
More precisely Gritsenko has shown that:
\bigskip
${\A}_{d}$ is not unirational 
(in fact $p_g(\tilde{\A}_{d})\ge 1$) if $d\ge 13$ and 
$d\ne 14,\; 15,\; 16,\; 18,\; 20,\; 24,\; 30,$ $36$. 
\bigskip

Moreover, a combination of the results in [Gri1] and [HS]
shows that ${\A}^{lev}_{p}$ is of general type for 
all primes $p\ge 37$.

In this paper, we will focus on some of the moduli spaces 
excluded from Gritsenko's analysis. In particular, we define certain moduli
spaces $\A_d^H$ of abelian surfaces with {\it partial level structure},
determined by groups $H$, which fit in between $\A^{lev}_d$ and 
$\A_d$: there are forgetful maps
$$\A_d^{lev}\rightarrow \A_d^{H}\rightarrow\A_d.$$
We will give details as to the structure of $\A_d$ (or rather
$\A_{d}^{H}$ for certain choices of $H$ depending on the case), 
for $d=12, 14, 16, 18$ and $20$.  In particular,
we will prove the following:

\theorem{results} For suitable choices of $H$, depending on the case,
\item{a)} $\A_{12}^{lev}$ is birational to a rational singular
complete intersection of two quadrics in $\Pfive$. 
\item{b)} $\A_{14}^{H}$ is birational to 
the complete intersection of the standard Pl\"ucker quadric and the
secant variety of the Veronese surface in $\P^5$. 
In particular it is unirational, being  a conic bundle with a rational
2-section. It is not rational.
\item{c)} $\A_{16}^{H}$ is birational to
$\Pone\times\Pone\times\Pone$, and thus rational.
\item{d)} $\A_{18}^{H}$ is birational
to a rational singular complete intersection of 2 quadrics in $\P^5$.
\item{e)} $\A_{20}^{H}$ is
rational.

\medskip
The basic method of proof for $d=12,14,16$ and $18$ is as follows.
For $d$ even, let
$Z_{d}=\bigcup_A A\cap \P_-$, where the union runs over all Heisenberg
invariant abelian surfaces of degree $(1,d)$ in $\P^{d-1}$, and
$\P_-$ is the $-1$ eigenspace of the standard involution $\iota:\P^{d-1}
\rightarrow\P^{d-1}$ inducing negation on $A$. For each such $A$,
$A\cap \P_-$ consists of the four odd two-torsion points
of $A$, and these four points form an orbit under the $\boldz_2\times\boldz_2$
action on $\P_-$ induced by $\sigma^{d/2}$ and $\tau^{d/2}$,
where $\sigma$ and $\tau$ are the standard generators of the Heisenberg
group $\HHH_{d}$. It was proved in
[GP1] that each such orbit in $Z_{d}$  is, in general, the
set of odd 2-torsion points of a unique Heisenberg invariant
abelian surface $A$, so that
$Z_{d}/\boldz_2\times\boldz_2$ is birational to $\A_{d}^{lev}$.
In the $d=12$ case, this quotient can be described explicitly. In higher
degrees, it is unlikely that this quotient is of negative Kodaira dimension, 
and in each of the higher degrees $14,16$ and $18$, some additional
particular geometry of the situation is used to take a further
quotient of $Z_{d}/\boldz_2\times\boldz_2$ to obtain a variety
birational to $\A_{d}^{H}$. The meaning of this moduli space is
described in \S 1.

The case $d=20$ is treated differently. Here we use the fact that there
exists singular quintics in $\Pfour$ which contain pencils of
(non-minimal) $(1,20)$-polarized abelian surfaces. This yields a ruling
of $\A_{20}^H$.

Of course such quintics are Calabi-Yau threefolds, and there exist small
resolutions of these quintics giving a non-singular Calabi-Yau with
Hodge number $h^{1,1}=2$ and $h^{1,2}=2$. This adds to the stable
of examples of Calabi-Yau threefolds containing a pencil of abelian
surfaces found in [GP3]. These example have recently been of interest
to physicists looking for interesting compactifications of space-time,
see the work of Candelas and Davies [CD] and Bak, Bouchard and Donagi
[BBD].

In fact, this Calabi-Yau threefold was the first
of the examples discovered in this project, found
by the second author in [Pop] while working on
the classification of surfaces in $\Pfour$. We realised that the existence of
this Calabi-Yau proved uniruledness of $\A_{20}$, so we embarked on
a search for other such Calabi-Yau threefolds, which led to the project
of which this is the final part.

There is one other example of a Calabi-Yau threefold arising in this paper,
a linear section of $\Gr(2,7)$, containing a pencil of $(1,14)$
abelian surfaces. In fact it turns out to be birational to
the Pfaffian Calabi-Yau of degree 14 in $\Psix$ containing a pencil
of $(1,7)$ abelian surfaces. There is some beautiful geometry associated
to this Calabi-Yau which we discuss in \S 3.
\bigskip

{\it Acknowledgments:}  It is a pleasure to thank Kristian Ranestad,
who joined one of us in discussions leading to some of the ideas
in this paper, and Ciro Ciliberto, Wolfram Decker, Igor Dolgachev,
David Eisenbud, Klaus Hulek, Nicolae Manolache, Gregory Sankaran, 
Frank Schreyer
and Allesandro Verra from whose ideas the exposition has benefited.
We are also grateful to Mike Stillman, Dave Bayer, and Dan Grayson for
the programs {\it Macaulay\/} [BS], and {\it Macaulay2\/} [GS]
which have been extremely useful to us; without them we
would perhaps have never been bold enough to guess the
existence of the structures that we describe here. Finally, we thank
Philip Candelas, without whose continued interest in Calabi-Yau manifolds
with small Hodge numbers, this paper would probably have remained
uncompleted forever.

\section {prelim} {Preliminaries.}

We review our notation and conventions concerning abelian surfaces; more
details can be found in [GP1].

Let $(A,\L)$ be a general abelian surface with a polarization of type
$(1,d)$.  If $d\ge 5$, then $|\L|$ induces an embedding of
$A\subset\P^{d-1}=\P(\dual{H^0(\L)})$ of degree $2d$.  The line bundle
$\L$ induces a natural map from $A$ to its dual $\hat A$,
$\phi_{\L}:A\longrightarrow \hat A$, given by $x\mapsto
t_x^*\L\otimes\L^{-1}$, where $t_x:A\longrightarrow A$ is the morphism
given by translation by $x\in A$.  Its kernel $K(\L)$ is isomorphic to
$\boldz_{d}\times\boldz_{d}$, and is dependent only on the
polarization $c_1(\L)$.

For every $x\in K(\L)$ there is an isomorphism $t_x^*\L
\cong \L$. This induces a projective representation
$K(\L)\longrightarrow {\PGL}(H^0(\L))$, which lifts uniquely to a
linear representation after taking a central extension
of $K(\L)$
$$1\longrightarrow{{\CC}^*}\longrightarrow{\G(\L)}
\longrightarrow{K(\L)}\longrightarrow0,
$$
whose Schur commutator map is the Weil pairing.  $\G(\L)$ is the {\it
theta group} of $\L$ and is isomorphic to the abstract Heisenberg
group $\H(1,d)$, while the above linear representation is isomorphic to the
Schr\"odinger representation of $\H(1,d)$ on $V={\CC}(\boldz_{d})$,
the vector space of complex-valued functions on $\boldz_{d}$.  An
isomorphism between $\G(\L)$ and $\H(1,d)$ which restricts to the
identity on centers induces a symplectic isomorphism between $K(\L)$
and $K(1,d)=\boldz_{d}\times\boldz_{d}$. This is an isomorphism which takes
the Weil pairing on $K(\L)$ to the standard skew-symmetric pairing $e^{(1,d)}$ on
$K(1,d)$ with values in $\CC^*$ given by
$$e^{(1,d)}\big((1,0),(0,1)\big)=\exp(-2\pi i/d).
$$
Such an isomorphism is called a {\it
level structure of canonical type} on $(A,c_1(\L))$. (See [LB],
Chapter 8, \S 3 or [GP1], \S 1.)

A decomposition $K(\L)=K_1(\L)\oplus K_2(\L)$, with $K_1(\L)\cong
K_2(\L) \cong \boldz_{d}$ isotropic subgroups with respect to the Weil
pairing, and a choice of a characteristic $c$ ([LB], Chapter 3, \S 1)
for $\L$, define a unique natural basis $\setdef{\vt^c_x}{ x\in K_1(\L)}$ of
{\it canonical theta functions} for the space $H^0(\L)$ (see [Mu2]
and [LB], Chapter 3, \S 2). This basis allows an identification of
$H^0(\L)$ with $V$ via $\vt_{\gamma}^c\mapsto x_{\gamma}$, where
$x_{\gamma}$ is the function on $\boldz_{d}$ defined by
$x_{\gamma}(\delta)=\cases{1&$\gamma=\delta$\cr
0&$\gamma\not=\delta$\cr}$ for $\gamma,\delta\in\boldz_{d}$. The functions
$x_0,\ldots,x_{d-1}$ can also be identified with coordinates on
$\P(\dual{H^0(\L)})$. Under this identification, the representation
$\G(\L)\longrightarrow \GL(H^0(\L))$ coincides with the Schr\"odinger
representation $\H(1,d)\longrightarrow \GL(V)$. We will only consider the
action of $\HHH_{d}$, the finite subgroup of $\H(1,d)$
generated in the Schr\"odinger representation by $\sigma$ and $\tau$,
where
$$ \sigma(x_i)=x_{i-1}, \qquad \tau(x_i)={\xi}^{-i} x_i,$$ for all
$i\in{\boldz}_{d}$, and where $\xi=e^{2{\pi}i\over d}$ is a primitive
root of unity of order $d$. Notice that $\lbrack\sigma,\tau\rbrack=\xi^{-1}$,
so ${\HHH}_{d}$ is a central extension
$$1 \longrightarrow {\bf\mu_{d}} \longrightarrow {\HHH_{d}}
\longrightarrow
{\boldz_{d}}\times{\boldz_{d}} \longrightarrow 0.$$
Therefore the choice of a canonical level structure means that
if $A$ is embedded in $\P(\dual{H^0(\L)})$ using as coordinates
$x_{\gamma}=\vt_{\gamma}^c$,  $\gamma\in
\boldz_{d}$, then the image of $A$ will be invariant under the
action of the Heisenberg group $\HHH_{d}$ via the Schr\"odinger
representation. (See [LB], Chapter 6, \S 7 for details.)

If moreover the line bundle $\L$ is chosen to be symmetric (and there
are always finitely many choices of such an $\L$ for a given
polarization type), then the embedding via $|\L|$ is also invariant
under the involution $\iota$, where
$$ \iota(x_i)=x_{-i}, \qquad i\in{\boldz}_{d}.$$
This involution restricts to $A$ as the involution $x\mapsto -x$.
We will denote  by $\P_+$ and $\P_-$ the $(+1)$ and
$(-1)$-eigenspaces of the involution $\iota$, respectively.
We will also write as usual
$\HHH_d^e:=\HHH_d\rtimes\langle\iota\rangle$.

We also recall a key result from [GP1]: In that paper, on $\P^{2d-1}\times
\P^{2d-1}$, we have introduced a matrix
$$M_d(x,y)=(x_{i+j}y_{i-j}+x_{i+j+d}y_{i-j+d})_{0\le i,j\le d-1},$$
where the indices of the variables $x$ and $y$ above are all modulo
$2d$. This matrix has the property that if $A\subseteq \P^{2d-1}$
is a Heisenberg invariant $(1,2d)$-polarized abelian surface,
then $M_d$ has rank at most two on $A\times A\subseteq \P^{2d-1}\times
\P^{2d-1}$. Similarly, if $A\subseteq \P^{2d}$
is a Heisenberg invariant $(1,2d+1)$-polarized abelian surface,
then the (Moore) matrix
$$M_{2d+1}'(x,y)=(x_{d(i+j)}y_{d(i-j)})_{i\in\boldz_{2d+1},j\in\boldz_{2d+1}},$$
on $\P^{2d}\times \P^{2d}$  has rank at most four
on $A\times A\subseteq \P^{2d}\times\P^{2d}$.
These matrices will prove to be ubiquitous!
\bigskip

{\bf Partial Level Structures.}

Recall from above 
that a canonical level structure on a polarized abelian surface
$(A,\L)$ with polarization of type $D=(d_1,d_2)$
is a symplectic isomorphism $b:K(\L)\rightarrow K(D)$.
(See also [GP1], \S 1.3 and [GP3] for notation and details.)

\definition{pls} Let  $H\subseteq K(D)$ be a fixed subgroup. An $H$-level 
structure on a $D$-polarized abelian surface $(A,\L)$ is an $H$-equivalence
class of canonical level structures $[b]$, where two canonical level
structures $b_1$ and $b_2$ are $H$-equivalent if there is a
symplectic automorphism $\phi:K(D)\rightarrow K(D)$ with $\phi|_H=id_H$ and
$\phi\circ b_1=b_2$.

The following proposition uses the notation of [GP1], \S 1.3.

\proposition{partiallev} Let $X_Z$ and $X_{Z'}$ be two fibres 
of $\X_D\rightarrow \H_g$, and suppose 
$\L|_{X_Z}$ and $\L|_{X_Z'}$ are both very ample,
so we can identify $X_Z$ and $X_{Z'}$ with their images in $\P(V)$, where
$V=H^0(\L|_{X_Z})=H^0(\L|_{X_{Z'}})$. Let $H\subseteq K(D)$ be a subgroup, and let
$H'$ be the inverse image of $H$ in $\H(D)$. We will think of $H'$ as
$H'\subseteq \GL(V)$ via the Schr\"odinger representation of $\H(D)$. 
Then the canonical level structures on $X_Z$ and $X_{Z'}$ are
$H$-equivalent if and only if there exists an element $T\in N(\H(D))$, the
normalizer of $\H(D)$ in $\GL(V)$, such
that $[T,\alpha]=T\alpha T^{-1}\alpha^{-1}\in {\bf C}^*$ for all $\alpha\in
H'$, and $T(X_Z)=X_{Z'}$. 

\proof: Suppose that the level structures $b$, $b'$ on $X_Z$ and $X_{Z'}$
are $H$-equivalent. Then there exists a symplectic automorphism $\phi
\in Sp(D)$ with $\phi\circ b=b'$ and $\phi|_H=id$. We have an exact 
sequence
$$0\rightarrow \H(D)\rightarrow N(\H(D))\rightarrow Sp(D)\rightarrow 1$$
where $N(\H(D))$ is the normalizer of $\H(D)$ in $\GL(V)$. Lift $\phi$ to
an element $T\in N(\H(D))$; then $T(X_Z)=X_{Z'}$, by
[GP1], Proposition 1.3.1. Also, the action of $\phi$ on $K(D)$ is given by
conjugation by $T$, so for any $\alpha\in H'$, $T\alpha T^{-1}$ must be 
a multiple of $\alpha$, so 
$T\alpha T^{-1}\alpha^{-1}\in {\bf C}^*$
for all $\alpha\in H'$.

Conversely, if there exists a transformation $T\in N(\H(D))$ with 
$T\alpha T^{-1}\alpha^{-1}\in {\bf C}^*$
for all $\alpha\in H'$, and $T(X_Z)=X_{Z'}$, then $X_Z$ and $X_{Z'}$ are
isomorphic as polarized abelian surfaces. 
$T$ then induces an element of $Sp(D)$
which is the identity on $H$, so that the canonical level structures on
$X_Z$ and $X_{Z'}$ are $H$-equivalent. \Box

With notation as in \ref{partiallev} we now introduce the following:

\definition{glh} If $H\subseteq K(D)$ is a subgroup, and
$\iota:V\rightarrow V$ the usual involution, 
define $\GL_H\subseteq \GL(V)$
by 
$$\GL_H=\{T\in \GL(V)\mid\hbox{$T\alpha T^{-1}\alpha^{-1}\in {\bf C}^*$
for all $\alpha \in H'$ and $\iota T=T \iota$}\},$$
where $H'$ is as in \ref{partiallev}.

In two special cases relevant for us we can explicitly describe
$\GL_H$:
 
\proposition{exampleglh} Suppose $D=(1,2d)$.
\item{\rm (1)}
If $H\subseteq K(D)$ is
the subgroup generated by $\sigma^2$ and $\tau$, then $T\in \GL_H$
if and only if either $T=\diag(a,b,a,b,\ldots,a,b)$ or $T=\sigma^d\circ
\diag(a,b,a,b,\ldots,a,b)$ for some $a,b\in{\bf C}^*$.
\item{\rm (2)} If $H=2K(D)$, then $T=(a_{ij})\in  \GL(V)$ is
in $\GL_H$ if and only if $a_{ij}=0$ unless $i=j$ or $i=j+d$, and
$a_{i+2,j+2}=a_{i,j}$. 
\item{\rm (3)} In either of the two above cases, if
$T\in \GL(V)$ satisfies $T \alpha T^{-1}\alpha^{-1}\in {\bf C}^*$
for all $\alpha\in H'$, then there exists $a,b$ such that, with
$T'=\sigma^a\tau^b T$, we have $T'\in \GL_H$.

\proof: Note that if $T=(a_{ij})$, and $T\circ \sigma^2=C_1
\sigma^2\circ T$ for some $C_1\in {\bf C}^*$, then 
$$a_{i-2,j-2}=C_1a_{ij}\quad\forall i,j\in\boldz/2d\boldz.\leqno{(1.1)}$$
Note that assuming $a_{ij}\not=0$ for some $i,j$, using this $d$ times
implies $C_1^d=1$, i.e., $C_1=\xi^{2n}$ for some $n\in \boldz/2d\boldz$,
where $\xi$ is a fixed primitive $2d$th root of unity.
If $T \circ \tau=C_2 \tau\circ T$, then
$$\xi^{-j}a_{ij}=C_2\xi^{-i}a_{ij}\quad \forall i,j\in\boldz/2d\boldz,
\leqno{(1.2)}$$
while if $T\circ \tau^2=C_3\tau^2\circ T$, then 
$$\xi^{-2j}a_{ij}=C_3\xi^{-2i}a_{ij}\quad \forall i,j\in\boldz/2d\boldz.
\leqno{(1.3)}$$
Finally, if $T\circ\iota=\iota\circ T$, then
$$a_{ij}=a_{2d-i,2d-j}.\leqno{(1.4)}$$

If (1.1) and (1.2) hold for $T$, then we see from (1.2) that 
either $a_{ij}=0$ or
$C_2=\xi^{i-j}$. Thus $C_2$ must be of the form $C_2=\xi^m$ for some
$m\in \boldz/2d\boldz$, and then
$a_{ij}=0$ unless $i-j=m$. 

If furthermore $T\in \GL_H$ for $H$ generated by $\sigma^2$ and $\tau$,
then (1.4) implies
$2m\equiv 0\mod 2d$, so $m=0$ or $d$, in which case $C_2=1$ or $-1$.
If $m=0$, then (1.1) and (1.4) applied to $(i,j)=(d+1,d+1)$ imply
together that $C_1=1$ and $T=\diag(a,b,a,b,\ldots)$. If $m=d$, then
(1.1) and (1.4) applied to $(i,j)=(1,d+1)$ imply together that $C_1=1$ and
$T=\sigma^d\circ \diag(a,b,a,b,\ldots)$.
This proves (1).

For (3) in the case that $H$ is generated by $\sigma^2$ and $\tau$,
note, from the properties of $T$ derived from (1.1) and (1.2), 
that if we set
$T'=\sigma^m\tau^{-n} T$,
and if $\{x_i\,|\,i\in \boldz/2d\boldz\}$ is the standard basis of $V$ as usual,
then
$$\eqalign{
T'(x_i)&=\sigma^m\tau^{-n}(a_{i+m,i}x_{i+m})\cr
&=\xi^{n(i+m)}a_{i+m,i}\sigma^m(x_{i+m})\cr
&=\xi^{n(i+m)}a_{i+m,i} x_i\cr
}$$
so the matrix $(a'_{ij})$ for $T'$ is diagonal, and using (1.1) and
$C_1=\xi^{2n}$,
$$a'_{ii}=\xi^{n(i+m)}a_{i+m,i}=\xi^{n(i+m)}\xi^{-2n}a_{i+m-2,i-2}
=\xi^{n(i+m)}\xi^{-2n}\xi^{-n(i+m-2)}a'_{i-2,i-2}=a'_{i-2,i-2}.$$
Thus $T'\in\GL_H$ by (1).

For (2), note that from (1.3), either
$a_{ij}=0$ or $C_3=\xi^{2i-2j}$. Thus $C_3$ must be of the form
$C_3=\xi^{2m}$, and then $a_{ij}=0$ unless $2(i-j)\equiv 2m\mod 2d$. 
But if $T\in \GL_H$, then (1.4) implies $4m\equiv 0\mod 2d$. We thus need to 
consider two cases:

{\it $d$ odd:} $2m=0$, and $a_{ij}=0$ unless $i-j=0$ or $d$. Applying
(1.4) to $(i,j)=(d+1,d+1)$ and $(1,d+1)$ shows that $C_1=1$ and we have the
desired form for $T$.

{\it $d$ even:} We have two cases: $2m=0$ or $2m=d$. In the first case,
it follows as in the odd case that $T$ is of the desired form, and we need
to rule out $2m=d$. In this case, $a_{ij}=0$ unless $i-j=d/2$ or
$3d/2$. But note by (1.4) that $a_{0,d/2}=a_{0,3d/2}$ and
$a_{d,d/2}=a_{d,3d/2}$. These are the only possible non-zero
entries in columns $d/2$ and $3d/2$, and since these columns are
then identical, $T$ cannot be invertible. Thus this case doesn't arise.

Finally, we show (3) in the case that $H$ is generated by $\sigma^2$
and $\tau^2$, so we assume $T$ is given satisfying (1.1) and (1.3).
Again setting $T'=\sigma^m\tau^{-n} T$, we see that
$$\eqalign{T'(x_i)&=\sigma^m\tau^{-n} (a_{i+m,i}x_{i+m}+a_{i+m+d,i}x_{i+m+d})\cr
&=\xi^{n(i+m)}a_{i+m,i}x_i+\xi^{n(i+m+d)}a_{i+m+d,i}x_{i+d}.\cr
}$$
Thus if $T'=(a'_{ij})$, then $a'_{ij}=0$ unless $i=j$ or $i=j+d$.
Furthermore, using (1.1) and $C_1=\xi^{2n}$, with $\delta=0$ or $d$,
$$\eqalign{a'_{i+\delta,i}&=\xi^{n(i+m+\delta)}a_{i+m+\delta,i}
=
\xi^{n(i+m+\delta)}\xi^{-2n}a_{i+m+\delta-2,i-2}\cr
&=\xi^{n(i+m+\delta)}\xi^{-2n}\xi^{-n(i+m+\delta-2)}a'_{i+\delta-2,i-2}
=a'_{i+\delta-2,i-2}.\cr}
$$
This shows from (2) that $T'\in\GL_H$.
\Box

Continuing to focus on the $(1,2d)$ case, recall from [GP1], \S 6
the set $Z_{2d}\subseteq \P_-$, which is the union of odd 2-torsion
points of Heisenberg invariant abelian surfaces of type $D=(1,2d)$. 
$\overline{Z_{2d}}$ is an irreducible $3$-fold for $d\ge 5$,
as was proven in [GP1], \S 6.

\definition{partiallevmod}
$\A_{2d}^H$ denotes the moduli space of $(1,2d)$-polarized abelian
surfaces with an $H$-level structure, for $H\subseteq K(D)$ a subgroup.
In particular, if $H=K(D)$, $\A^H_{2d}=\A^{lev}_{2d}$, the moduli
space of $(1,2d)$-polarized abelian surfaces with canonical level structure.

The reason for introducing
$\GL_H$ is the following theorem:

\theorem{partialmod} If $d\ge 5$ and $H\subseteq K(D)$ is either
generated by $\sigma^2$ and $\tau$ or is $2K(D)$, then $\A^H_{2d}$
is birationally equivalent to
$\overline{Z_{2d}}/(\GL_H\cap N(\H(1,2d)))$.

\proof: 
We know by [GP1], Theorem 6.2 that there is a birational
map $\Theta_{2d}:\A_{2d}^{lev}
\rto\overline{Z_{2d}}/\boldz_2\times\boldz_2$ taking $[A]
\in\A^{lev}_{2d}$ to the $\boldz_2\times\boldz_2=\langle\sigma^d,\tau^d
\rangle$-orbit $A\cap\P_-$, the odd two-torsion points of $A$. One notes
from \ref{exampleglh} that $\sigma^d,\tau^d\in GL_H\cap N(\H(1,2d))$, so that
we get a factorization of the quotient 
$$\overline{Z_{2d}}\rightarrow\overline{Z_{2d}}/\boldz_2\times\boldz_2
\rightarrow\overline{Z_{2d}}/(GL_H\cap N(\H(1,2d))).$$
In particular, the action of $GL_H\cap N(\H(1,2d))$ descends to
$\overline{Z_{2d}}/\boldz_2\times\boldz_2$.

Letting $X_Z$, $X_{Z'}$ be two general fibres of the universal family
$\X_{(1,2d)}\rightarrow
\H_2$ as in \ref{partiallev}, these represent the same point in $\A^H_{2d}$
if and only if $T(X_Z)=X_{Z'}$ for some $T\in N(\H(1,2d))$ such that
$T\alpha T^{-1}\alpha^{-1}\in{\bf C}^*$ for all $\alpha\in H'$. Thus, if
there exists $T\in GL_H\cap N(\H(1,2d))$ such that $T(\Theta_{2d}([X_Z]))
=\Theta_{2d}([X_{Z'}])$, it follows that $T(X_Z\cap\P_-)=X_{Z'}\cap\P_-$.
Indeed, $T$ commutes with $\iota$, and hence takes odd two-torsion points
of $X_Z$ to odd two-torsion points of $X_{Z'}$.
Since $T$ is in the normalizer of $\H(1,2d)$, $T(X_Z)$ is Heisenberg
invariant, and by [GP1], Theorem 6.2, the general
Heisenberg invariant abelian surface is determined by its intersection with
$\P_-$. Thus $T(X_Z)=X_{Z'}$ and $X_Z$ is $H$-equivalent to $X_{Z'}$.

Conversely, if $X_Z$ is $H$-equivalent to $X_{Z'}$, we obtain 
$T\in N(\H(1,2d))$ with $T\alpha T^{-1}\alpha^{-1}\in {\bf C}^*$ for all 
$\alpha\in H'$, but we need not have $T\iota=\iota T$. By \ref{exampleglh},
we have $T'=\sigma^a\tau^b T$ which does satisfy $T'\iota=\iota T'$.
So $T'\in
GL_H\cap N(\H(1,2d))$ and $T'([X_Z])=[X_{Z'}]$. Thus $\overline{Z_{2d}}
/(GL_H\cap N(\H(1,2d)))$ is birationally equivalent to $\A^H_{2d}$, as desired.
\Box

\section{1.12} {Moduli of $(1,12)$-polarized abelian surfaces.}

Given the results of [GP1], \S 2 and \S 6, it is not difficult to
describe the structure of $\A_{12}^{lev}$. By [GP1], Corollary 2.7, if
$A\subseteq\P^{11}$ is an $\HHH_{12}$-invariant abelian surface with
polarization of type $(1,12)$, $y\in A$ a point, then $M_6(x,y)$ has
rank at most two on $A$, so in particular, if $x\in A\cap \Pfour_-$,
the matrix $M_6(x,x)$ must have rank at most two. Now if
$x\in\Pfour_-$, $M_6(x,x)$ is a $6\times 6$ skew-symmetric matrix, and
thus $x\in A\cap\Pfour_-$ for some Heisenberg invariant abelian
surface implies that the $4\times 4$ Pfaffians of $M_6(x,x)$ vanish.

For $x\in \Pfour_-$ with coordinates $$
\hbox{$(x_1:x_2:x_3:x_4:x_5)$ ($=
(0:x_1:x_2:x_3:x_4:x_5:0:-x_5:-x_4:-x_3:-x_2:-x_1)\in \P^{11}$)},$$
the first 
$4\times 4$ block of $M_6(x,x)$ is
$$\pmatrix{
0&-x_1^2-x_5^2&-x_2^2-x_4^2&-2x_3^2\cr
x_1^2+x_5^2&0&-x_1x_3-x_3x_5&-2x_2x_4\cr
x_2^2+x_4^2&x_1x_3+x_3x_5&0&-2x_1x_5\cr
2x_3^2&2x_2x_4&2x_1x_5&0\cr}$$
and the $4\times 4$ pfaffian of this is
$$f=2x_3^2(x_1x_3+x_3x_5)-2(x_2^2+x_4^2)x_2x_4+2x_1x_5(x_1^2+x_5^2).$$
Let $Q$ be the quartic hypersurface in $\Pfour_-$ determined by the equation
$f=0$.

\theorem {rat.map12} 
The map $\Theta_{12}:\A_{12}^{lev}\rto \Pfour_-/\boldz_2
\times\boldz_2$ induces a birational map to $Q/\boldz_2\times\boldz_2$.

\proof: By [GP1], Theorem 6.2, $\Theta_{12}$ is birational onto its image.
Let $\pi:\Pfour_-\rightarrow \Pfour_-/\boldz_2\times\boldz_2$ be the 
projection, and let $Z_{12}$ be the inverse image under $\pi$ of $\im
\Theta_{12}$. Then by construction, $\overline{Z_{12}}\subseteq Q$. On the
other hand, $\overline{Z_{12}}$ must be three dimensional. It is easy to check
that $Q$ is a non-singular quartic hypersurface, and thus $\overline{Z_{12}}=Q$.
\Box

We now give a more precise description of $Q/\boldz_2\times\boldz_2$.

\theorem{descrip12} $Q/\boldz_2\times\boldz_2$ is isomorphic
to the complete intersection in $\Pfive$ given by the equations
$$\eqalign{z_0z_5-z_1z_4+z_2z_3&=0 \quad\quad\hbox{(the Pl\"ucker quadric)}\cr
z_0z_2-z_3^2+2z_2z_5&=0\cr}$$
In particular, $Q/\boldz_2\times\boldz_2$ is rational.
\par

\proof: We first note that the $\boldz_2\times\boldz_2$ action on
$\Pfour_-$ takes the form
$$\eqalign{\sigma^6:(x_1:x_2:x_3:x_4:x_5)&\mapsto (x_5:x_4:x_3:x_2:x_1)\cr
\tau^6:(x_1:x_2:x_3:x_4:x_5)&\mapsto (x_1:-x_2:x_3:-x_4:x_5)\cr}$$
and the ring of invariants of ${\bf C}[x_1,\ldots,x_5]$ under this
action is generated by $x_3,x_1+x_5,x_1^2+x_5^2,x_2^2+x_4^2,$ and $x_2x_4$,
and thus
the invariant quadrics
$$\eqalign{z_0=x_1^2+x_5^2&\quad z_3=x_1x_3+x_3x_5\cr
z_1=x_2^2+x_4^2&\quad z_4=x_2x_4\cr
z_2=x_3^2&\quad z_5=x_1x_5\cr}$$
generate the invariants of even degree of this action. Thus
$\Pfour_-/\boldz_2\times\boldz_2$ is isomorphic to a hypersurface in
$\Pfive=\Proj {\bf C}[z_0,\ldots,z_5]$. It is easy to see that the quadrics
$z_0,\ldots,z_5$ satisfy the relation
$$z_0z_2-z_3^2+2z_2z_5=0.$$
On the other hand, $Q$ is given by 
$$f=2(z_0z_5-z_1z_4+z_2z_3).$$
This gives the desired equations.

Rationality of this threefold can be seen easily by
adapting the standard proof of rationality (see e.g. [GH],
pg. 796) of the non-singular quadric line complex to this singular
one. One projects from a line contained in $Q/\boldz_2\times\boldz_2$
to $\Pthree$. For example, one can use the line given by 
$z_0=z_1=z_3=z_5=0$.
\Box

\section {1.14} {Moduli of $(1,14)$-polarized abelian surfaces.}

This is the first case considered here
in which the locus of odd 2-torsion points
(in $\P_-$) is no longer of negative Kodaira dimension. Following the
strategy for $(1,12)$, one finds that the matrix $M_7(x,x)$ has rank
2 on a quartic hypersurface $f_1=0$ in $\Pfive_-$, and the matrix
$M_7(\sigma(x),x)$ has rank 2 on another quartic hypersurface $f_2=0$.
One finds $f_1=f_2=0$ is an irreducible threefold, and hence is the closure
of the locus of odd 2-torsion points, $\overline{Z_{14}}$ in
the notation of [GP1], \S 6. The quotient of this threefold
by $\boldz_2\times\boldz_2$ is birational to $\A^{lev}_{14}$, and is
very likely of general type, but it is difficult to determine this. So at this
point, we need to stop considering the moduli space with full level 
structure.
Instead, we will consider the subgroup $H\subseteq K(1,14)$ generated
by $\sigma^2$ and $\tau$, and find the moduli space $\A^H_{14}$ of
$(1,14)$ abelian surfaces with an $H$-level structure. 

We will first study the geometry of certain mappings which will enable us
to divide out $\overline{Z_{14}}$ by $\boldz_2\times \boldz_2$ as well as make 
the additional identifications which identify $H$-equivalent level structures.
It will turn out that this is easier than computing $\overline{Z_{14}}/
\boldz_2\times \boldz_2$. 

More precisely, 
looking at the matrix $M_7(x,x)$ on $\Pfive_-$, one sees that it 
contains a submatrix
$$M'=\pmatrix{
0 & -x_1^2-x_6^2 & -x_2^2-x_5^2 & -x_3^2-x_4^2\cr
x_1^2+x_6^2& 0 & -x_1x_3-x_4x_6 & -x_2x_4-x_3x_5\cr
x_2^2+x_5^2 & x_1x_3+x_4x_6 & 0 & -x_1x_5-x_2x_6\cr
x_3^2+x_4^2 & x_2x_4+x_3x_5 & x_1x_5+x_2x_6&0\cr}$$
which is rank $\le 2$ on a quartic $Q_1=\{f_1=0\}$ where $f_1$
is the $4\times 4$ Pfaffian of $M'$. Similarly, $M_7(\sigma(x),x)$ 
restricted to $\Pfive_-$ has a submatrix
$$M'_{\sigma}=\pmatrix{
0& -x_1x_2-x_5x_6 & -x_2x_3-x_4x_5 & -2x_3x_4\cr
x_1x_2+x_5x_6& 0  & -x_1x_4-x_3x_6 & -2x_2x_5\cr
x_2x_3+x_4x_5& x_1x_4+x_3x_6&0& -2x_1x_6\cr
2x_3x_4& 2x_2x_5 & 2x_1x_6 & 0 \cr}$$
$M'_{\sigma}$ has rank 2 on the set
$Q_2=\{f_2=0\}$, where $f_2$ is (one-half of) the Pfaffian of $M'_{\sigma}$. 
Since $M'_{\sigma}$ has
rank 2 at the general point of $Q_2$, we obtain a rational map
$\phi:Q_2\rto \Gr(2,4)$ induced by $M'_{\sigma}$, taking a point
of $Q_2$ to the kernel of $M'_{\sigma}$ at that point.
This is defined where $M'_{\sigma}$ 
is non-zero,
and in Pl\"ucker coordinates is given by
$$(x_1,\ldots,x_6)\mapsto (z_0,\ldots,z_5)=
(
x_1x_2+x_5x_6,
x_2x_3+x_4x_5, x_1x_4+x_3x_6,
2x_3x_4, 2x_2x_5,  2x_1x_6).$$

Here,
$$f_2=x_1x_3x_4^2-x_2^2x_3x_5-x_2x_4x_5^2+x_1^2x_2x_6+x_3^2x_4x_6
    +x_1x_5x_6^2$$ 
(this is in fact the unique bi-degree $(2,2)$ bi-invariant of $\PSL_2(\boldz_7)$
acting on the variables $x_1,x_3,x_5$ and $x_2,x_4,x_6$ diagonally in
its three-dimensional representation)
and
$$f_1-f_2=x_1x_3^3-x_2^3x_4+x_1^3x_5-x_3x_5^3+x_4^3x_6+x_2x_6^3,$$
(which is in fact the sum of two Klein quartics, one in the variables $x_1,x_3,
x_5$, the other in the variables $x_2,x_4,x_6$, the Klein quartic being the
unique quartic invariant of $\PSL_2(\boldz_7)$ in its three-dimensional
representation).

\lemma{geom14} \item{\rm (1)} $\phi$ is not defined on the union
of the planes $P_1=\{x_1=x_3=x_5=0\}$ and $P_2=\{x_2=x_4=x_6=0\}$.
\item{\rm (2)} The closure of the  image of $\phi$,
$\overline{\im\phi}$, is a complete intersection of type $(2,3)$
whose equations are the Pl\"ucker quadric $z_0z_5-z_1z_4+z_2z_3=0$ and
the cubic 
$$\det\pmatrix{ z_5&z_2&z_0\cr
z_2&z_3&z_1\cr
z_0&z_1&z_4\cr}=0.$$
This latter hypersurface is the secant variety of the Veronese surface
determined by the $2\times 2$ minors of this matrix.
\item{\rm (3)} $\GL_H$ acts on the fibres of $\phi$,
and the fibre of $\phi$ over a general point is the closure
of a $\GL_H$-orbit.
\item{\rm (4)} $Q_2$ is irreducible.
\item{\rm (5)} $\overline{\im \phi}$ is an irreducible, unirational $3$-fold.

\proof: (1) One checks easily that the zero locus of the six quadrics defining
the map $\phi$ is $P_1\cup P_2$.

(2) Consider first the rational map $\phi':\Pfive_-\rto \P^8$ induced
by the linear system $|H^0(\I_{P_1\cup P_2}(2))|$. Specifically,
using coordinates $\{z_{ij}\}$ on $\P^8$ with $i\in\{1,3,5\}$ and
$j\in \{2,4,6\}$, we give $\phi'$ via $z_{ij}=x_ix_j$. The equations
of $\overline{\im\phi'}$ are then given as the $2\times 2$ minors of the matrix
$$L=\pmatrix{ z_{16} & z_{14} & z_{12}\cr
z_{36} & z_{34} & z_{32}\cr
z_{56} & z_{54} & z_{52}\cr}.$$
If $X$ is the blow-up of $\Pfive_-$ along $P_1\cup P_2$, then $\phi'$
lifts to a morphism $\tilde\phi':X\rightarrow\P^8$ which is easily
seen to describe $X$ as a $\Pone$-bundle over $\overline {\im\phi'}$,
which is naturally isomorphic to $P_1\times P_2$ embedded via the
Segre embedding. Explicitly, a line $l$ in $\Pfive_-$ meeting $P_1$ and
$P_2$ gets mapped to the point $(p_1,p_2)$ of $P_1\times P_2$
with $l\cap P_1=\{p_1\}$ and $l\cap P_2=\{p_2\}$. 

Now $\phi$ is the composition of $\phi'|_{Q_2}$ and the linear projection
$\pi:\P^8\rto\Pfive$ given by $\pi:(z_{ij})\mapsto (z_0,\ldots,z_5)
=(z_{12}+z_{56}, z_{32}+z_{54},z_{14}+z_{36}, 2z_{34}, 2z_{52}, 2z_{16})$.
The center of the projection $\pi$ is easily seen to be disjoint from
$\overline{\im\phi'}$, so $\pi$ must map $\overline{\im\phi'}$ to
a four-fold in $\Pfive$. In fact, this projection identifies points
on $P_1\times P_2\subseteq \P^8$ via the involution which exchanges
$P_1$ and $P_2$, $P_1$ and $P_2$ being identified via $\sigma^7$ acting
on $\Pfive_-$. One sees this by noting that the equations for $\pi$
are invariant under the induced action of $\sigma^7$ on $\P^8$.
Now since $L$ has rank 1 on 
$\overline{\im\phi'}$, 
$$L+L^t=\pmatrix{ 2z_{16} & z_{14}+z_{36} & z_{12}+z_{56}\cr
z_{14}+z_{36} & 2z_{34} & z_{32}+z_{54}\cr
z_{12}+z_{56} & z_{32}+z_{54} & 2z_{52}\cr}$$
has rank $\le 2$ on $\overline{\im\phi'}$. Thus
$$\det\pmatrix{ z_5&z_2&z_0\cr
z_2&z_3&z_1\cr
z_0&z_1&z_4\cr}=0$$
is an equation for $\overline{\im\phi}$. As is well-known, see for instance 
[SR] or [GH], pages 179--180,
this determines an irreducible cubic hypersurface which is the secant variety
of the Veronese surface determined by the $2\times 2$ minors of
this matrix. Thus 
$\overline{\im \pi\circ\phi'}$ must coincide with this cubic hypersurface,
which is then isomorphic to $P_1\times P_2/\boldz_2$.
Since $Q_2$ can be expressed in terms of $z_0,\ldots,z_5$ as the Pl\"ucker
quadric, $\pi\circ \phi'(Q_2)$ then is the complete intersection of this
Pl\"ucker quadric and the cubic hypersurface, showing (2).

To complete the proof of (3), note that $\pi:\overline{\im\phi'}
\rightarrow \overline{\im\pi\circ\phi'}$ must be a 2:1 map. Thus the fibres
of $\pi\circ\phi'$ consist of pairs of lines. But in addition, one
easily sees from the equations of $\phi$ and the structure of $\GL_H$
from \ref{exampleglh} that the action of $\GL_H$ on $\Pfive_-$ restricts
to an action on the fibres of $\phi$, from which (3) follows. 

(4) To prove that $Q_2$ is irreducible, first note that the inverse
image of $\overline{\phi(Q_2)}\subseteq P_1\times P_2/\boldz_2$ in 
$P_1\times P_2$ is given by the intersection of the quadric
$$2z_{16}(z_{12}+z_{56})-2z_{52}(z_{32}+z_{54})+2z_{34}(z_{14}+z_{36})=0$$
and $P_1\times P_2\subseteq \P^8$. This coincides with $\overline
{\phi'(Q_2)}$. This divisor of type $(2,2)$ in $P_1\times P_2$ is easily seen
to be irreducible by inspection. Now the
proper transform of $Q_2$ in $X$ is a $\Pone$-bundle over 
$\overline{\phi'(Q_2)}$, so $Q_2$ must be irreducible.

(5) By the above discussion, 
$\overline{\im\phi}$ is covered by $\overline{\phi'(Q_2)}$,
which is an irreducible divisor of type $(2,2)$ in $P_1\times P_2$.
Thus $\overline{\im\phi}$ is irreducible. It is easy to see that
an irreducible divisor $D$ in $P_1\times P_2$ of type $(2,2)$ is unirational:
$D$ is a conic bundle via $p_2:P_1\times P_2\rightarrow P_2$,
while if $l\subseteq P_1$ is a general line, then $p_1^{-1}(l)\cap D$
is a rational 2-section of $p_2$, and thus $D$ is unirational:
see e.g. [Be], Cor.\ 4.4. \Box

\theorem{A14moduli} $\overline{\im\phi}$ is birational to $\A^H_{14}$.

\proof: We know that $\overline{Z_{14}}\subseteq \{f_1=f_2=0\}\subseteq
\Pfive_-$. (We have not proven that the latter algebraic set is irreducible,
so we don't know if this is an equality.) Let $x\in\overline{Z_{14}}$
be a general point, and consider the orbit $\GL_H\cdot x\subseteq
\Pfive_-$ of $x$. Since $x\in Q_2$, by \ref{geom14} (3), $\GL_H\cdot x\subseteq
Q_2$. Note that for general $x$, $\overline{\GL_H\cdot x}$ consists of two
lines, say $l_1^x\cup l_2^x$. Where does $l_1^x$ intersect $\overline{
Z_{14}}$? Since $\tau^7\in \GL_H$, both $x$ and $\tau^7(x)\in l_1^x$. 
In addition, if $T=\diag(i,1,i,1,\ldots)$, one sees that
$T\in N(\H(1,14))$, so $T(x)$ and $\tau^7(T(x))$ are both in
$l_1^x\cap{\overline{Z_{14}}}$. Furthermore, for general $x$, these
four points are distinct. Similarly, $l_2^x=\sigma^7(l_1^x)$ intersects $
\overline{Z_{14}}$ in at least four distinct points. If $l_1^x$ intersected
$\overline{Z_{14}}$ in additional points, so would $l_2^x$, and vice versa,
and then we would necessarily have $\overline{\GL_X\cdot x}\subseteq
\{f_1=f_2=0\}$, which is a threefold as $Q_2$ is irreducible, by \ref{geom14},
(4). Note that if $\overline{\GL_H\cdot x}\subseteq \{f_1=f_2=0\}$, then
in particular, the fact that $f_1-f_2$ vanishes on the point
$(ax_1,bx_2,ax_3,bx_4,ax_5,bx_6)$ for all $a$ and $b$ implies that
$k_1:=x_1x_3^3+x_1^3x_5-x_3x_5^3=0$
and $k_2:=-x_2^3x_4+x_4^3x_6+x_6^3x_2=0$, i.e., $x$ is in the join
of the Klein quartic curves in $P_1$ and $P_2$. Clearly
$Q_2$ does not contain this join, and hence $x$ is contained in 
a surface $Q_2\cap\{k_1=k_2=0\}$, so $x$ is not a general point of
$\overline{Z_{14}}$.
Thus for general $x$, $l_1^x$ and $l_2^x$
intersect $\overline{Z_{14}}$ in precisely four points.
In particular, $\overline{\GL_H\cdot x}\cap \overline{Z_{14}}$
is precisely the orbit of $x$ under the group $\GL_H\cap N(\H(1,14))$. 

Thus, by \ref{geom14} (3), $\overline{\phi(Z_{14})}$ is birational
to $\overline{Z_{14}}/(\GL_H\cap N(\H(1,14)))$, which in turn is
birational to $\A_{14}^H$ by \ref{partialmod}. But since $\overline{Z_{14}}
\subseteq Q_2$ and $\overline{\phi(Z_{14})}\subseteq \overline{\im\phi}$,
and both of these latter varieties are irreducible threefolds, we must
have equality, proving the theorem. \Box 

\remark{irr14} A more careful analysis of $\overline{\im\phi}$ in 
fact shows that this threefold has a resolution of singularities
which is a conic bundle over $\Ptwo$, with discriminant locus a non-singular
plane
curve of degree 6. It then follows from [Be], Theorem 4.9, that this threefold
is irrational. Hence unirationality is the best we can achieve here.

\definition{V14y} Let $y\in\P^5_-$. Define $V_{14,y}\subseteq \P^{13}$
to be the scheme defined by the $4\times 4$ Pfaffians of $M_7(x,y)$.

Recall also the definition of the Pfaffian Calabi-Yau threefold in
$\Psix$ from [GP3], \S 5:

\definition{V7y} Let 
$$M_7'(x,y)=(x_{3(i+j)}y_{3(i-j)})_{i,j\in \boldz_7}$$
on $\Psix\times\Psix$. For $y\in\Ptwo_-\subseteq\Psix$, we define 
$V_{7,y}$ to be the scheme defined by the $6\times 6$ Pfaffians of
$M_7'(x,y)$.

In [GP3], \S 5, $V_{7,y}$ was seen to be, for general $y$, a degree 14
Calabi-Yau threefold with 49 ordinary double points, the latter being the
$\HHH_7$-orbit of $y$.

$V_{14,y}$ is a linear section of $\Gr(2,7)$, and there is a correspondence
between Pfaffian Calabi-Yau threefolds in $\Psix$ and such linear
sections of $\Gr(2,7)$ in $\P^{13}$ which now appears to be
folklore, see for example [R\o d]. The general
construction is as follows.

Let $V$ be a 7-dimensional vector space over ${\bf C}$, and let
$W\subseteq \bigwedge^2 V$ be a seven-dimensional subspace. Corresponding
to this subspace is a subscheme $Pf_W$ given by the intersection
$\P(W)\cap D_4\subseteq \P(\bigwedge^2 V)$, where $D_4$ is the projectivized
locus of rank $\le 4$ elements of $\bigwedge^2 V$. This scheme $Pf_W$ is,
for general choice of $W$, a Pfaffian Calabi-Yau threefold of degree 14 in 
$\P(W)$.

On the other hand, the annihilator $W^0\subseteq \bigwedge^2\dual{V}$ of
$W$ is a 14 dimensional subspace. Let $Gr_{W^0}=\P(W^0)
\cap D_2^*\subseteq\P(\bigwedge^2\dual{V})$ where $D_2^*$ is the locus
of rank $\le 2$ alternating forms on $V$. Of course $D_2^*\cong \Gr(2,7)$
via the Pl\"ucker embedding. For general $W$, $Gr_{W^0}$ is a non-singular
degree 42 Calabi-Yau threefold in $\P^{13}$. 

\remark{} It has been observed in [R\o d] that these two families of
Calabi-Yau threefolds have the same mirror. This led us to
conjecture in 1997 that $Pf_W$ and $Gr_{W^0}$ have isomorphic
derived categories of coherent sheaves. This has now been proved by
Borisov and C\u{a}ld\u{a}raru [BC] and Kuznetsov [K]. 
See also [HT] for a physics argument, and
the delightful term ``glop'' for this correspondence.

In what follows, $\Pfive_-$ now denotes the negative eigenspace
of $\iota$ acting on $\P^{13}$, and $\Ptwo_-$ denotes that of
$\iota$ acting on $\Psix$.

\proposition{corr714} Let $\psi:\Pfive_-\rto\Ptwo_-$ be defined by taking
$x=(x_1,\ldots,x_6)$ to the kernel of the matrix
$$\pmatrix{x_5&x_3&x_1\cr x_2&x_4&x_6\cr},$$
i.e., 
$$\psi(x_1,\ldots,x_6)=(x_3x_6-x_1x_4,x_1x_2-x_5x_6,x_4x_5-x_2x_3).$$
Then
\item{(1)} The fibres of $\psi$ are $\Pthree$'s which are closures of
$GL^{2K(1,14)}$-orbits.
\item{(2)} Under the correspondence between Pfaffian and Grassmann
Calabi-Yau threefolds, $V_{14,y}$ corresponds to $V_{7,\psi(y)}$.

\proof: (1) follows easily from the description of $GL^{2K(1,14)}$
in \ref{exampleglh}.

(2) Let $y'=(y_0',\ldots,y_6')\in \Ptwo_-\subseteq\Psix=\P(V)$,
with $y'_i=-y'_{-i}$. Then $V_{7,y'}$ is given by the ideal generated
by $6\times 6$ Pfaffians of $M_7'(x,y')$. Let
$z_{i,j}$ be coordinates on $\P^{20}=\P(\bigwedge^2 V)$, with
$0\le i,j\le 6$, $z_{i,j}=-z_{j,i}$, so that $(z_{i,j})_{0\le i,j\le 6}$
is the corresponding skew-symmetric matrix. Then $M'_7(x,y')$ defines
a map $\alpha:\Psix\rightarrow \P^{20}$ via
$$z_{i,j}=x_{3(i+j)}y'_{3(i-j)}.$$
Now $\alpha$ is, for general choice of $y'$, an embedding, and $V_{7,y'}
=\alpha^{-1}(D_4)$.

Now we wish to compute the equations of the annihilator of
$\alpha(\Psix)$ in $\P(\bigwedge^2 \dual{V})$. Note that $\alpha(\Psix)$
is spanned by the images of $e_0,\ldots,e_6$, (the standard basis
vectors), and $e_i$ is mapped via $\alpha$ to the point with
coordinates
$$z_{j,k}=\delta_{i,3(j+k)}y'_{3(j-k)}.$$
The ideal of the annihilator is generated by the hyperplanes in
$\P(\bigwedge^2\dual{V})$ corresponding to these seven points. If
$z_{i,j}^*$ are the dual coordinates on $\P(\bigwedge^2\dual{V})$,
then this gives us the ideal for the annihilator generated by the seven
equations
$$y_1'z_{1+i,6+i}^*+y_2'z_{2+i,5+i}^*+y_3'z_{3+i,4+i}^*=0,\quad 0\le i\le 6.$$
This defines $L\subseteq\P(\bigwedge^2 \dual{V})$, $L\cong\P^{13}$, and
$L\cap D_2^*$ is the Grassmann Calabi-Yau corresponding to $V_{7,y'}$.
To compare this with $V_{14,y}$ for some $y$, we have to find
a parametrization $\beta:\P^{13}\rightarrow L$ given by the
matrix $$M_7(x,y)=(x_{i+j}y_{i-j}+x_{i+j+7}y_{i-j+7})_{0\le i,j\le 6}.$$
In other words, a point $y\in \Pfive_-$ yields a linear map
$\beta:\P^{13}\rightarrow\P(\bigwedge^2\dual{V})$ by taking the point
$x=(x_i)\in \P^{13}$ to the point with coordinates $(z_{i,j}^*)$ given
by
$$z_{i,j}^*=x_{i+j}y_{i-j}+x_{i+j+7}y_{i-j+7},\quad 0\le i,j\le 6,$$
and $V_{14,y}=Gr_{\im\beta}$. In order for the image of $\beta$ to be $L$,
the seven equations defining $L$ must be satisfied by any point in
the image of $\beta$, i.e.,
$$y_1'(x_{7+2i}y_9+x_{2i}y_2)+y_2'(x_{7+2i}y_{11}+x_{2i}y_4)
+y_3'(x_{7+2i}y_{13}+x_{2i}y_6)=0$$ 
for $0\le i\le 6$. Since this must hold for all values of $x$, and keeping
in mind in addition that $y\in\Pfive_-$, we find we must have
$$\eqalign{y_1'y_5+y_2'y_3+y_3'y_1&=0,\cr
y_1'y_2+y_2'y_4+y_3'y_6&=0.\cr}$$
This is precisely equivalent to the statement
$$(y_1',y_2',y_3')\in \ker\pmatrix{y_5&y_3&y_1\cr y_2&y_4&y_6\cr}.$$
In other words, $y'=\psi(y)$. It is also easy to see that
the map $\beta$ defined by such a choice of $y$  is an isomorphism
with $L$ as long as $(y_5,y_3,y_1)$ and $(y_2,y_4,y_6)$ are linearly
independent. \Box

\remark{V14yremarks} Using Macaulay or Macaulay 2 [BS], [GS], one can check
the following three facts:
\item{(1)} For general $y\in\Pfive_-$, $V_{14,y}$
contains a pencil of $(1,14)$-polarized abelian surfaces invariant
under the subgroup $\langle\sigma^2,\tau^2\rangle$ of $\HHH_{14}$.
\item{(2)} For general $y\in\Pfive_-$, $V_{14,y}$ is an irreducible
threefold of degree 42 with 49 ordinary double points.
\item{(3)} If $y\in Z_{14}$ is general, then the $6\times 6$
Pfaffians of $M_7(x,\sigma^7(y))$ induce a birational map
$V_{14,y}\rto V_{7,\psi(y)}$.

Unfortunately, we do not know a non-computational proof of (2) and (3). 
It is intriguing
that the correspondence between $Pf_W$ and $Gr_{W^0}$ identifies, in
the particular case of $V_{14,y}$ and $V_{7,\psi(y)}$, birational
Calabi-Yau threefolds. This is not true in general: it is easy to
see that for general $W$, $Pf_W$ and $Gr_{W^0}$ are not birationally
equivalent, as they both have Picard number one but for $H$ a primitive
generator of the Picard group, $H^3=14$ and $42$ respectively.

Note that (3) implies that the Calabi-Yau $V_{14,y}$, or equivalently,
$V_{7,\psi(y)}$, is birationally fibred by abelian surfaces in two different
ways, one fibration having $(1,7)$-polarized fibres, and one having
$(1,14)$-polarized fibres. The existence of these two distinct pencils
of abelian surfaces was observed in Remark 2.9 of [MR].

\section {1.16} {Moduli of $(1,16)$-polarized abelian surfaces.}

Restricting $M_8(x,x)$ to $\P^6_-$, one finds a $5\times 5$ block
$$M=\pmatrix{
0  & -x_1^2-x_7^2 & -x_2^2-x_6^2 & -x_3^2-x_5^2 & -2x_4^2\cr
x_1^2+x_7^2 & 0 & -x_1x_3-x_5x_7 & -x_2x_4-x_4x_6 & -2x_3x_5\cr
x_2^2+x_6^2 & x_1x_3+x_5x_7 & 0 & -x_1x_5-x_3x_7 & -2x_2x_6\cr
x_3^2+x_5^2 & x_2x_4+x_4x_6 & x_1x_5+x_3x_7  &  0 & -2x_1x_7\cr
2x_4^2 & 2x_3x_5 & 2x_2x_6 & 2x_1x_7 & 0\cr}$$
In the discussion that follows, it will be helpful to make a change
of coordinates. We define a new set of coordinates $(y_0,y_1,z_0,z_1,
t_0,t_1,u)$ on $\Psix_-$ by
$$y_0=x_1+x_7, y_1=x_3+x_5, z_0=x_1-x_7, z_1=x_3-x_5, t_0=x_4, t_1=x_2+x_6,
u=x_2-x_6.$$
In these coordinates, 
$$M={1\over 2}\pmatrix{
0&-y_0^2-z_0^2&-t_1^2-u^2&-y_1^2-z_1^2&-4t_0^2\cr
y_0^2+z_0^2&0&-y_0y_1-z_0z_1&-2t_0t_1&-y_1^2+z_1^2\cr
t_1^2+u^2&y_0y_1+z_0z_1&0&-y_0y_1+z_0z_1&-t_1^2+u^2\cr
y_1^2+z_1^2&2t_0t_1&y_0y_1-z_0z_1&0&-y_0^2+z_0^2\cr
4t_0^2&y_1^2-z_1^2&t_1^2-u^2&y_0^2-z_0^2&0\cr}$$

Let $Z\subseteq \P^6_-$ be the subvariety of defined by the $4\times 4$
Pfaffians of $M$. Clearly $\overline{Z_{16}}\subseteq Z$.

\lemma{deg40} $Z$ is an irreducible threefold of degree $40$ in $\P^6_-$.
In particular $\overline{Z_{16}}=Z$.

\proof: If $Z$ is in fact a threefold, then $Z$ is of the expected
degree for such a Pfaffian variety, and standard degree
calculations show it must be of degree $40$. To show it is an irreducible 
threefold, we must work harder. First, consider the point
$$p=(y_0,y_1,z_0,z_1,t_0,t_1,u)=(0,0,0,0,0,0,1)\in\Psix_-.$$
The tangent cone to $Z$ at $p$ can be computed by working on the
affine open set of $\Psix_-$ where $u=1$, and identifying the leading terms
of the $4\times 4$ Pfaffians of $M$ with $u$ replaced by $1$. These leading
terms yield the following equations for the tangent cone at $p$:
$$2t_0t_1=y_0^2-z_0^2+y_1^2+z_1^2=y_0^4-y_1^4-z_0^4+z_1^4+8t_0^3t_1
=y_0^2+y_1^2+z_0^2-z_1^2=0.$$
Note however that $y_0^4-y_1^4-z_0^4+z_1^4+8t_0^3t_1$ is in the ideal
generated by the three quadrics, so in fact the tangent cone has
ideal 
$$I=(t_0t_1,y_0^2+y_1^2,z_0^2-z_1^2).$$
This is in fact a union of 8 affine 3-planes. 

Now consider the projection map $\pi:\Psix_-\rto \Pfive$ given
by $(y_0,y_1,z_0,z_1,t_0,t_1,u)\mapsto (y_0,y_1,z_0,z_1,t_0,t_1)$
and its restriction $\pi:Z\rto \Pfive$. First note that $\pi$
does not map any component of $Z$ to something of lower dimension.
If it did, the component would have to be a cone with vertex at $p$.
Such a component would then appear as a component of the tangent
cone. However, one checks easily that none of the 8 three-planes in the
tangent cone at $p$ are in fact contained in $Z$. To understand
the image $\pi(Z)$, we first find two equations vanishing on $\pi(Z)$.
Note that the Pfaffian of the submatrix of $M$ obtained by deleting
the third row and third column of $M$ does not depend on
$u$, and this Pfaffian vanishes on $\pi(Z)$. It is
$$g_1=y_0^4-y_1^4-z_0^4+z_1^4+8t_0^3t_1.$$
For another equation, one adds the Pfaffian obtained by deleting the
first row and column of $M$ with the Pfaffian obtained by deleting
the fifth row and fifth column, to get
$$g_2=2y_0^3y_1+2y_0y_1^3-2z_0^3z_1+2z_0z_1^3-4t_0t_1^3.$$
Thus $\pi(Z)\subseteq Z'=\{g_1=g_2=0\}$. Now
the singular locus of $Z'$ can be computed by hand,
and one finds in fact that $Z'$ is a codimension 2 complete intersection
(hence of degree 16) with 6 singular points, at $(1,\pm i,0,0,0,0)$,
$(0,0,1,\pm 1,0,0)$, $(0,0,0,0,1,0)$, and $(0,0,0,0,0,1)$. Thus in 
particular $Z'$ is irreducible, so $\pi(Z)=Z'$. We also see that 
$\pi$ maps $Z$ two-to-one to $Z'$, to account for the degree of
$Z'$: $\deg Z'=(\deg Z-\mult_P Z)/\deg\pi$. In fact, $Z$ is clearly
invariant under the involution $j:\Psix_-\rightarrow\Psix_-$ given by negation
of the coordinate $u$, and this interchanges the two sheets of the double
covering $\pi:Z\rto Z'$. Thus $Z$ consists of either one or two
irreducible components. If it is the latter, then $j$ interchanges
these two components. Note that in coordinates $x_1,\ldots,x_7$,
$j$ acts by interchanging $x_2$ and $x_6$.

To show that $Z$ is irreducible, we note that if not, then the branch
locus $B\subseteq Z'$ of the double covering $p:Z\rto Z'$
must be non-reduced everywhere. This branch locus is the image of
$Fix(j)\cap Z$, which is computed easily by setting $u=0$ in
the matrix $M$ and taking $4\times 4$ Pfaffians to obtain the equations
for $B\subseteq Z'$. We then just need to identify one point
in $B$ where the tangent cone to $B$ is reduced. Such a point
is $(y_0,\ldots,t_1)=(0,0,0,0,0,1)$, which is contained in $B$
and whose tangent cone is easily seen, as we did earlier, to be
given by 
$$t_0=y_0^2-y_1^2=z_0^2+z_1^2=0,$$
which is reduced. \Box

\remark{generaltype} With a bit of additional calculation, one finds the
6 singularities of $Z'\subseteq\Pfive$ are ordinary triple points (i.e.
have tangent cone a cone over a non-singular cubic surface).
In particular these are canonical singularities, Thus $Z'$ is in fact
of general type with canonical bundle $\O_{Z'}(2)$. Thus $Z=\overline{Z_{16}}$
is also of general type.

Our goal now is to understand $\A_{16}^H$ with $H=2K(1,16)$. In this case,
$\GL_H$ is a four-dimensional group, as computed in \ref{exampleglh}.
First, we will understand the set of three-dimensional $\GL_H$ orbits
in $\P^6_-$. 

For $x=(x_1,\ldots,x_7)\in \P^6_-$, let $N_x$ be the matrix
$$N_x=\pmatrix{-z_1 & 0 & z_0 & 0 & -z_0 & 0 & z_1\cr
0 & t_0 & 0 &-t_1 & 0 & t_0 & 0\cr
y_1 & 0 & -y_0 & 0 & -y_0 & 0 & y_1\cr}.$$
We can then define a rational map
$\phi:\P^6_-\rto \Gr(4,7)$ by $x\mapsto \ker N_x$.
Here the columns of $N_x$ correspond to the coordinates $x_1,\ldots,x_7$.

\lemma{geom116} 
\item{\rm (1)} $rank(N_x)<3$ if and only if $x\in \Pi_1\cup\Pi_2\cup\Pi_3$
where
$$\eqalign{\Pi_1&=\{z_0=z_1=0\},\cr
\Pi_2&=\{t_0=t_1=0\},\cr
\Pi_3&=\{y_0=y_1=0\}.\cr}$$
\item{\rm (2)} $\dim \overline{\GL_H\cdot x}=3$ if and only if $x\not\in\Pi_1
\cup\Pi_2\cup\Pi_3\cup \{u=0\}$.
\item{\rm (3)} $\overline{\GL_H\cdot x}=\P(\ker N_x)$ if $x\not\in\Pi_1\cup\Pi_2
\cup\Pi_3\cup\{u=0\}$.
\item{\rm (4)} $\overline{\im\phi}$ is birational
to $\Pone\times\Pone\times\Pone$.

\proof: (1) Clearly $N_x$ drops rank on $\Pi_1\cup\Pi_2\cup\Pi_3$. Conversely,
if $N_x$ is rank $\le 2$, then either the second row of $N_x$ is
zero or the first and third rows are linearly dependent. This latter case
occurs only on $\Pi_1\cup\Pi_3$.

(2) 
By \ref{exampleglh}, one sees that
the action of $\alpha\in \GL_H$ on $\P^6_-$ is of the form
$$\eqalign{
\alpha:(x_1,\ldots,x_7)\mapsto&(ax_1-bx_7,cx_2-dx_6,\cr
&ax_3-bx_5,(c-d)x_4,ax_5-bx_3,\cr
&cx_6-dx_2,ax_7-bx_1)\cr}$$
for suitable $a,b,c,d\in {\bf C}$. In our new coordinates this is
$$\alpha:(y_0,y_1,z_0,z_1,t_0,t_1,u)\mapsto ((a-b)y_0,(a-b)y_1,
(a+b)z_0,(a+b)z_1,(c-d)t_0,(c-d)t_1,(c+d)u).$$
Thus $\GL_H\cong ({\bf C}^*)^4$, and $\dim\overline{\GL_H\cdot x}<3$
if and only if $x$ is in $\Pi_1$, $\Pi_2$, $\Pi_3$, or $\{u=0\}$.
(2) then follows.

(3) We first note that $N_x\cdot x=0$, so $x$ is in the $\Pthree$
represented by $\phi(x)$. 
From the form of $\alpha$ given in (2),
one sees that $N_{\alpha(x)}$ is obtained from $N_x$ by multiplying
the first row of $N_x$ by $a+b$, the second row by $c-d$ and the third by
$a-b$. Thus as $N_{\alpha(x)}\cdot\alpha(x)=0$, we also have
$N_x\cdot\alpha(x)=0$, so $\overline{\GL_H\cdot x}\subseteq \P(\ker N_x)$.
By (1) and (2), these are the same dimension so equality holds.

(4) Let $L_1,L_2,L_3\subseteq \dual{\P^6_-}$ be lines dual to
$\Pi_1$, $\Pi_2$ and $\Pi_3$ respectively. If $p\in\im\phi$, then $p$ 
corresponds to a $\Pthree\subseteq \P^6_-$, with dual $\Ptwo$ being,
say, $P\subseteq \dual{\P^6_-}$. Then $P\cap L_i$  consists
of one point for each $i$. This gives a birational map between
$\overline{\im\phi}$ and $L_1\times L_2\times L_3$. \Box 

\lemma{number116} $$\#\left({\GL_H\cap N(\H(1,16))\over {\bf C}^*}\right)=32.$$

\proof: 
First, $GL_H\cap \H(1,16)$ is generated by ${\bf C}^*$, $\sigma^8$ and $\tau^8$,
as is easily seen from \ref{exampleglh}. Taking the quotient of $N(\H(1,16))$ by
$\H(1,16)$ gives $SL_2(\boldz_{16})$.
Thus it is enough to identify the image of $GL_H\cap N(\H(1,16))$ in
$SL_2(\boldz_{16})$, i.e., as in the proof of \ref{partiallev}, we just need
to identify the subgroup of $SL_2(\boldz_{16})$ keeping $2K(1,16)$ invariant.
These are the elements $\alpha$ of $SL_2(\boldz_{16})$, $\alpha=(\alpha_{ij})_{1
\le i,j\le 2}$, with $\alpha (2,0)=(2,0)$, $\alpha(0,2)=(0,2)$. Thus
$\alpha_{11}$ and $\alpha_{22}$ are either $1$ or $9$, and $\alpha_{12}$
and $\alpha_{21}$ are either $0$ or $8$. This gives a total of $16$ possibilities,
but of these, half have determinant $9$, leaving a total of $8$ such elements
with determinant $1$.
Since $\# (GL_H\cap\H(1,16))/{\bf C}^*=4$, this gives a total of $4\times 8=32$,
as desired.
\Box

\theorem{A116} $\overline{\im\phi}$ is birational to $\A_{16}^H$.

\proof: \ref{geom116} shows that $(\Psix_- \setminus(\Pi_1\cup\Pi_2\cup\Pi_3
\cup\{u=0\}))/\GL_H$ is birational to $\Pone\times\Pone\times\Pone$
via the map $\phi$. On the other hand, by \ref{partialmod},
$\A^H_{16}$ is birational to $Z_{16}/\GL_H\cap N(\H(1,16))$.
Thus to show $\A_{16}^H$ is birational to $\Pone\times\Pone\times\Pone$,
it is enough to show that for a general point $x\in Z_{16}$, 
$(\GL_H\cdot x)\cap Z_{16}=(\GL_H\cap N(\H(1,16)))\cdot x$. By
\ref{number116}, it is enough then to show that
$\#(\GL_H\cdot x)\cap Z_{16}=32$ for $x\in Z_{16}$ general. Since
$(\GL_H\cap N(\H(1,16)))\cdot x\subseteq (\GL_H\cdot x)\cap Z_{16}$,
clearly $\#(\GL_H\cdot x)\cap Z_{16}\ge 32$. On the other hand,
$\overline{Z_{16}}$ is degree 40, so if $\overline{\GL_H\cdot x}
\cap \overline{Z_{16}}$ is a finite set of points, we must then have
$$\#(\GL_H\cdot x)\cap Z_{16}\le 40.$$ Since this number must be divisible
by $32$, we see that $\#(\GL_H\cdot x)\cap Z_{16}=32$. Thus we are done
if we can find a point $x\in Z_{16}$ such that $\overline{\GL_H\cdot x}
\cap\overline{Z_{16}}$ is finite. Since the variety $\overline{Z_{16}}$ is
defined over the ring $\boldz[1/2]$, it is enough to find such a point
in finite odd characteristic. We know of no better way of finding such
a point except by searching exhaustively. For example, in characteristic
$p=7$, one finds the point with coordinates $y_0=4$, $y_1=-1$, $z_0=-2$,
$z_1=-1$, $t_0=2$, $t_1=-2$, $u=1$, does the trick. This can be checked
in Macaulay/Macaulay 2.
\Box

\section {1.18} {Moduli of $(1,18)$-polarized abelian surfaces.}

Again, we study $\A_{18}^H$ with $H=2K(1,18)$. As we shall see,
the closure of $\GL_H\cdot Z_{18}$ is a quartic hypersurface in 
$\P^7_-$, and we can use similar tricks to the $(1,14)$ case to understand
the correct quotient of $\overline{Z_{18}}$. 

The matrix $M_9(\sigma\tau^9(x),x)$ restricted to $\P^7_-$ has a $4\times
4$ block obtained by taking the $ij$-th entries with $1\le i\le 4$, $1\le j\le
4$, of the form
$$M':=\pmatrix{
0&x_1x_2-x_7x_8& -x_2x_3+x_6x_7& x_3x_4-x_5x_6\cr
-x_1x_2+x_7x_8&0 & x_1x_4-x_5x_8 & -x_2x_5+x_4x_7\cr
x_2x_3-x_6x_7& -x_1x_4+x_5x_8 &0& x_1x_6-x_3x_8\cr
-x_3x_4+x_5x_6 & x_2x_5-x_4x_7 & -x_1x_6+x_3x_8& 0\cr}.$$
Let $f$ be the Pfaffian of $M'$ and let $Q=\{f=0\}$. 
We also define a rational map $\phi:\P^7_-\rto \Gr(2,4)$
by taking a point $(x_1,\ldots,x_8)$ to the subspace of ${\bf C}^4$
spanned by the rows of the matrix
$$N=\pmatrix{x_1&x_3&x_5&x_7\cr
x_8&x_6&x_4&x_2\cr};
$$
explicitly in Pl\"ucker coordinates,
$$\eqalign{(x_1,\ldots,x_8)\mapsto (y_0,\ldots,y_5)=&
(x_1x_6-x_3x_8,-x_1x_4+x_5x_8,x_1x_2-x_7x_8,x_3x_4-x_5x_6,\cr 
&-x_2x_3+x_6x_7,x_2x_5-x_4x_7).\cr}$$

\lemma{geom18} \item{\rm (1)} $\phi$ fails to be defined precisely on the
fourfold $Z\subseteq \P^7_-$ defined by the $2\times 2$ minors of $N$;
this is isomorphic to $\Pone\times\Pthree$.
\item{\rm (2)} If $X$ is the blow-up of $\P^7_-$ along $Z$, then
$\phi$ lifts to a morphism $\tilde\phi:X\rightarrow \Gr(2,4)$
which is a surjection, a $\Pthree$ bundle, and each fibre of $\tilde\phi$
is the closure of a $\GL_H$ orbit.
\item{\rm (3)} If $\tilde Q$ is the proper transform of $Q$, then
$\tilde\phi(\tilde Q)$ is a (2,2) complete intersection defined
by the equations
$$\eqalign{y_0y_5-y_1y_4+y_2y_3&=0\cr
y_0y_2-y_1y_3+y_4y_5&=0\cr}$$
and $\tilde Q$ is a $\Pthree$-bundle over $\tilde\phi(\tilde Q)$
whose fibres are all closures of $\GL_H$-orbits.
\item{\rm (4)} $\tilde\phi(\tilde Q)$ is rational.

\proof: (1) is standard. For (2), $\phi$ clearly lifts to the
blowup, $\tilde\phi:X\rightarrow \Gr(2,4)$. It is a surjection
since any two-dimensional subspace of ${\bf C}^4$ can be expressed
as the row span of some $2\times 4$ matrix. Now an element $T=(a_{ij})
\in GL_H$ is determined by $\alpha=a_{00}$, $\beta=a_{11}$, $\gamma=a_{0d}$
and $\delta=a_{1,d+1}$, by \ref{exampleglh}, (2), and restricting the
action of $T$ to $\P^7_-$, one sees that the matrix $N$ is transformed by
$T$ to
$$\pmatrix{\beta x_1-\delta x_8&\beta x_3-\delta x_6&\beta x_5-\delta x_4&
\beta x_7-\delta x_2\cr
\alpha x_8-\gamma x_1&\alpha x_6-\gamma x_3&\alpha x_4-\gamma x_5&\alpha
x_2-\gamma x_7\cr}$$
which clearly has the same row span as $N$. So fibres of $\tilde\phi$
are $\GL_H$ invariant. In fact, two rank two $2\times 4$ matrices $N$
and $N'$ have the same row span if and only if $N'$ is obtained from 
$N$ in this fashion, for some choice of $\alpha,\beta,\gamma,\delta$,
necessarily with $\alpha\beta-\gamma\delta\not=0$. So the fibres
of $\tilde\phi$ are always closures of $\GL_H$ orbits. The closure
of the $\GL_H$ orbit of $(x_1,\ldots,x_8)$ is then the
span of the points $(x_1,0,x_3,0,x_5,0,x_7,0)$, $(x_8,0,x_6,0,x_4,0,x_2,0)$,
$(0,x_2,0,x_4,0,x_6,0,x_8)$ and $(0,x_7,0,x_5,0,x_3,0,x_1)$, hence
a $\P^3$ as long as this orbit is three-dimensional, which is
the case if $(x_1,\ldots,x_8)\not\in Z$. This shows (2).

(3) The first quadric is the Pl\"ucker relation defining $\Gr(2,4)$,
and the second equation just pulls back to the equation $f$. The
singularities of this complete intersection are easily computed
and found to consist of $6$ points; hence it must be
irreducible. Now the fibres of $\tilde\phi|_{\tilde Q}$ are
dimension $\le 3$, and $\tilde Q$ is dimension $6$, so in fact the
fibres of $\tilde\phi|_{\tilde Q}$ are three-dimensional
and $\tilde\phi(\tilde Q)$ is three-dimensional. Hence $\tilde\phi(\tilde Q)$
is precisely the complete intersection and the fibres are $\P^3$'s. 

(4) One checks easily that the point $P=(y_0,\ldots,y_5)=(1,0,0,1,1,0)$
is a double point; indeed, setting $y_0=1$ and expanding
the two equations in Taylor series shows the tangent cone to $\tilde\phi(\tilde
Q)$ at $P$ is a quadric cone. 
Consider the projection from $P$ given by
$$(y_0,\ldots,y_5)\mapsto (z_1,\ldots,z_5)=(y_1,y_2,y_3-y_0,y_4-y_0,y_5).$$
Noting that the difference of the two quadratic equations defining
$\tilde\phi(\tilde Q)$ is
$$-y_1(y_4-y_3)+y_2(y_3-y_0)+y_5(y_0-y_4)=-z_1(z_4-z_3)+z_2z_3-z_5z_4,$$
which defines a quadric cone $Q'$ in $\P^4$ with one singular point, we see
that the image of $\tilde\phi(\tilde Q)$ under this projection is contained
in $Q'$. Since 
$\tilde\phi(\tilde Q)$ is degree 4, there are two possibilities. Either
the projection $\tilde\phi(\tilde Q)\rto Q'$ is birational, or 
$\tilde\phi(\tilde Q)$ is a cone with vertex $P$. The latter is impossible,
since as observed earlier, $\tilde\phi(\tilde Q)$ has $6$ isolated
singular points.
\Box

We next need to refine the results of [GP1]:

\proposition{118eqs} For a general $\HHH_{18}$-invariant $(1,18)$-polarized
abelian surface $A\subseteq \P^{17}$, and a point $y\in A\cap
\P^7_-$, the $4\times 4$ Pfaffians of the matrix $M_9(\sigma\tau^9(x),y)$
generate the homogeneous ideal of $A$.

\proof: As in the proof of [GP1], Theorem 6.2, it is enough to show there is
a point $y\in\overline{Z_{18}}$ such that the $4\times 4$ Pfaffians of
the matrix $M(\sigma\tau^9(x),y)$ generate the homogeneous ideal of a surface
$S$ which has the same Hilbert polynomial as a $(1,18)$-polarized abelian
surface in $\P^{17}$. We proceed as follows to find such a point.

Let $E_1\subseteq\P^5$ be a general $\HHH_6$-invariant elliptic curve, and 
$E_2\subseteq \P^{17}$ be an $\HHH_{18}$-invariant elliptic curve,
thus $E_1$ and $E_2$
having canonical level structures of degrees 6 and 18 respectively.
Write coordinates $x_i$, $i\in\boldz_6$, on $\P^5$ and coordinates
$y_j,j\in\boldz_{18}$, on $\P^{17}$. Denote by $\sigma_6,\tau_6$ and
$\sigma_{18},\tau_{18}$ the generators of $\HHH_6$ and $\HHH_{18}$
respectively. Now $\boldz_6$ acts on $E_1\times E_2$
by $(x,y)\mapsto (\sigma^i_6(x),\sigma^{3i}_{18}(y))$ for $i\in\boldz_6$.
Then $A=(E_1\times E_2)/\boldz_6$ is an abelian surface which inherits
a polarization coming from the product polarization on $E_1\times E_2$.
Let $\pi:E_1\times E_2\rightarrow A$ be the quotient map.
The linear system of $\boldz_6$-invariant sections of the product
polarization on $E_1\times E_2$
induces a linear system on $A$, and this in fact gives
a (possibly rational) map $\alpha$ from $A$ into $\P^{17}(z)$
with coordinates $z_i,i\in\boldz_{18}$, by mapping
$E_1\times E_2\subseteq\P^5(x)\times\P^{17}(y)$ via
$$
z_i=\sum_{j\in\boldz_6} x_jy_{3j+i}.
$$
It is easy to see that these sections $z_i$ are $\boldz_6$-invariant.
Furthermore, the usual action of $\HHH_{18}$ on $\P^{17}(z)$
restricts to an action on $A$:
$$
\eqalign{
\sigma_{18}(\pi(x,y))&=\pi(x,\sigma_{18}(y))\cr
\tau_{18}(\pi(x,y))&=\pi(\tau_6^{-1}(x),\tau_{18}(y))\cr
}
$$
so clearly the quotient polarization $\L$ on $A$ 
is $\HHH_{18}$-invariant and hence comes with a
canonical level structure of type $(1,18)$. Now, in fact the linear system on 
$A$ given by $z_0,\ldots,z_{17}$ is base-point-free: by [LB], 10.1.1 and 10.1.2,
this linear system can only have base-points if there is an elliptic curve
$E\subseteq A$ with $E.c_1(\L)=1$. However such an elliptic curve must
then pull-back to a number of copies of elliptic curves on $E_1\times E_2$.
But the only elliptic curves on $E_1\times E_2$, for general non-isogenous
$E_1$, $E_2$, are of the form $E_1\times \{y\}$ or $\{x\}\times E_2$.
While $E_1\times \{y\}$ is degree 6 with respect to the product polarization,
its image in $A$ is also degree $6$, so in particular, there is no elliptic
curve $E\subseteq A$ with $E.c_1(\L)=1$. Furthermore, by Reider's theorem
([LB], 10.4.1), the morphism $\alpha:A\rightarrow \P^{17}$ is an embedding
as long as there is no elliptic curve $E\subseteq A$ with $E.c_1(\L)=2$;
this is indeed the case by a similar argument, so $\alpha$ is an embedding.

Next consider the 2-torsion points of $A$. Since $\boldz_6\subseteq E_1
\times E_2$ contains two two-torsion points, there are two sorts of
two-torsion points on $A$: First, there are $8$ points which are images of
$2$-torsion points on $E_1\times E_2$, and then there are 2-torsion points
which are images of points on $E_1\times E_2$ satisfying
$$-(x,y)=(\sigma_6(x),\sigma_{18}^3(y)).$$
(Here, the minus sign is in the group law on $E_1\times E_2$.)

In terms of projective coordinates, $x=(x_0,\ldots,x_5)$, $y=(y_0,\ldots,y_{17})$,
this is saying that $\iota(x)=\pm\sigma_6(x)$ and 
$\iota(y)=\pm\sigma_{18}^3(y)$, as the eigenvalues of $\sigma^{-1}_6\circ\iota$
and $\sigma^{-3}_{18}\circ\iota$ are $\pm 1$.
Note there are four solutions to the first equation on $E_1$,
and four solutions to the second on $E_2$, since an elliptic curve
has $4$ two-torsion points. If $\iota(x)
=\sigma_6(x)$, then $\iota(\tau^3_6(x))=-\sigma_6(\tau_6^3(x))$, and 
if $\iota(y)
=\sigma_{18}^3(y)$, then $\iota(\tau^9_{18}(y))=-\sigma^3_{18}
(\tau^9_{18}(y))$. Thus
solutions of these equations come in pairs, and there are two solutions each
to the four equations 
$$\iota(x)=\sigma_6(x), \quad \iota(x)=-\sigma_6(x),\quad
\iota(y)=\sigma_{18}^3(y),\quad \iota(y)=-\sigma^3_{18}(y).$$
Now if $\iota(x)=\sigma_6(x)$ and $\iota(y)=-\sigma^3_{18}(y)$, then
$z:=\alpha(x,y)$ satisfies $\iota(z)=-z$. Similarly, if $\iota(x)
=-\sigma_6(x)$, $\iota(y)=\sigma^3_{18}(y)$, then $z=\alpha(x,y)$
satisfies $\iota(z)=-z$. There are a total
of $8$ such pairs $(x,y)$, hence $4$ such points $z$, and we have
identified the odd $2$-torsion points of $A$. These points are in $Z_{18}$.

We now degenerate $A$ by degenerating the elliptic curve $E_{18}$ to
$$X(\Gamma_{18}):=\bigcup_{i\in\boldz_{18}} l_{i,i+1}\subseteq\P^{17},$$
where $l_{i,i+1}$ denotes the line with equations
$$
\{y_0=y_1=\cdots=y_{i-1}=y_{i+2}=\cdots=y_{17}=0\}
$$
(see [GP1], \S 4). It is well-known that such a curve is a degeneration
of $\HHH_{18}$-invariant elliptic curves. Thus if we now take
$A=(E_1\times X(\Gamma_{18}))/\boldz_6$, $A\cap\P^7_-\subseteq
\overline{Z_{18}}$. We thus just need to show that if $z'\in
A\cap\P^7_-$, then the $4\times 4$ Pfaffians of
$M_9(\sigma\tau^9(z),z')$ generate the homogeneous ideal of $A$.

Let us identify such a $z'$: this can be taken to be the
image of a point $(x,y)\in E_1\times E_2$, with $x\in E_1$ satisfying
$\iota(x)=\sigma_6(x)$ and $y\in E_2$ satisfying $\iota(y)
=-\sigma^3_{18}(y)$. We can write $x=(x_0,x_1,x_2,x_2,x_1,x_0)$, and
for $y=(y_i)$ we can take all coordinates $y_i=0$ except
$y_7=1,y_8=-1$. Then the image $z$ of $(x,y)$ is
$$
z=(0,x_2,-x_2,0,x_1,-x_1,0,x_0,-x_0,0,x_0,-x_0,0,x_1,-x_1,0,x_2,-x_2).
$$
We need to identify a choice of $x_0,x_1,x_2$ so that $x$ is on an 
$\HHH_6$-invariant elliptic curve. To do so, recall from [GP1], Corollary 2.2,
that the matrix $M_3(x,x)=(x_{i+j}x_{i-j}+x_{i+j+3}x_{i-j+3})_{0\le i,j\le 2}$
has rank $\le 1$ for any point $x\in \Pfive$ contained in a
$\HHH_6$-invariant elliptic normal curve. Restricting to the plane given
by $x_0=x_5, x_1=x_4, x_2=x_3$ gives the matrix
$$\pmatrix{x_0^2+x_2^2&(x_0+x_2)x_1&(x_0+x_2)x_1\cr
2x_1^2&2x_0x_2&2x_0x_2\cr
x_0^2+x_2^2&(x_0+x_2)x_1&(x_0+x_2)x_1\cr};
$$
there is one independent $2\times 2$ minor, namely 
$2x_0x_2(x_0^2+x_2^2)-2x_1^3(x_0+x_2)$. 

This defines a non-singular quartic in $\Ptwo$, and since the moduli
space of elliptic curves is one-dimensional, every point in this 
curve is a point on some such elliptic curve. Taking a specific
point on this curve, say in characteristic $p=11$, $(x_0,x_1,x_2)=(5,6,1)$,
yields a point $z'\in \overline{Z_{18}}$. One can then calculate
in Macaulay/Macaulay 2 that the Hilbert polynomial of the ideal
generated by the $4\times 4$ Pfaffians of the matrix $M_9(\sigma\tau^9(z),z')$
has the same Hilbert polynomials of a $(1,18)$-polarized abelian surface
in $\P^{17}$, as desired.

While the degenerate abelian surface described by this ideal is not difficult
to describe as a union of elliptic scrolls, it seems that we really
do need to resort to Macaulay to verify the statement about the ideal;
it seems there is no sufficiently simple choice of degeneration which
does the trick for us, for if $z\in\overline{Z_{18}}$ is chosen 
insufficiently generally, $M_9(\sigma\tau^9(z),z')$ does not 
generate an ideal with the right Hilbert polynomial. This explains
the difficulty
in proving this result without appeal to computational means.
\Box

\theorem{AH118theorem}
$\A^H_{18}$ is birational to $\tilde\phi(\tilde Q)$ (see \ref{geom18},
(3), (4)), and hence in particular $\A^H_{18}$ is rational.

\proof: \ref{geom18} shows that $(\tilde Q-Z)/\GL_H$ is birational
to $\tilde\phi(\tilde Q)$ via the map $\phi$. On the other hand,
by \ref{partialmod}, $\A^H_{18}$ is birational to
$Z_{18}/(\GL_H\cap N(\H(1,16)))$. Thus to show $\A^H_{18}$ is
birational to $\tilde\phi(\tilde Q)$, it is enough to show that for
a general point $y\in Z_{18}$, 
$$(\GL_H\cdot y)\cap Z_{18}=(\GL_H\cap N(\H(1,18)))\cdot y.$$
Suppose then that there is another point $y'\in (\GL_H\cdot y)
\cap Z_{18}$ such that $y'\not\in (\GL_H\cap N(\H(1,18)))\cdot y$.
Then write $y'=T(y)$ for some $T\in\GL_H$. Note by
\ref{exampleglh}, $T=(a_{ij})$ is determined by $\alpha=a_{0,0}$,
$\beta=a_{1,1}$, $\gamma=a_{0,d}$ and $\delta=a_{1,d+1}$.
It is then not hard to check explicitly that
$$
M_9(\sigma\tau^9(x),y')=M_9(\sigma\tau^9(x),T(y))=M_9(T'(\sigma\tau^9(x)),y),
$$
where $T'=(a'_{ij})\in\GL_H$ is determined by $a'_{0,0}=\beta$, 
$a'_{1,1}=\alpha$, $a'_{0,d}=-\gamma$, $a'_{1,d+1}=-\delta$. In fact,
$T'=(\alpha\beta-\gamma\delta)T^{-1}$. 

Since $y$ and $y'$ can be assumed to be sufficiently general, it follows
from \ref{118eqs} that the $4\times 4$ Pfaffians of $M_9(\sigma\tau^9(x),y)$
and $M_9(\sigma\tau^9(x),y')$ generate the ideals of $(1,18)$-polarized
abelian surfaces $A_y$ and $A_{y'}$ respectively. Thus if we view $T'$
as also inducing a map $T':\P^{17}\rightarrow\P^{17}$, we see that the 
$4\times 4$ Pfaffians of $M_9(T'(\sigma\tau^9(x)),y)$ generate
the ideal of $(T')^{-1}(A_y)=T(A_y)$. However, this ideal coincides
with that generating $A_{y'}$, so $A_{y'}=T(A_y)$, and $A_{y'}\cong
A_y$. This can only happen if $T\in N(\H(1,18))$ (see [GP1], Proposition 1.3.1),
contradicting
$y'\not\in (\GL_H\cap N(\H(1,18)))\cdot y$.
\Box

\section {1.20} {Moduli of $(1,20)$-polarized abelian surfaces.}

We deal with the (1,20) case using completely different methods,
relying largely on work carried out in the second author's thesis [Pop]
and in [ADHPR2]. This approach revolves around the fact that a certain 
blow-up of a general $(1,20)$-polarized abelian surface in 25 points
can be embedded in $\P^4$. The moduli space of these surfaces in
$\P^4$ is then seen to be rational.

We begin by reviewing some material from [ADHPR2]. In the discussion
which follows, we will have representations both of the group
$\HHH_{20}$ and the group $\HHH_5$. Denote the Schr\"odinger representations
of these two groups by $V_{20}$ and $V_5$ respectively.
We then have $\P^{19}=\P(\dual{V_{20}})$,
on which $\HHH_{20}$ acts, with eigenspaces of the standard involution
$\iota$ yielding $\P^8_-,\P^{10}_+\subseteq \P^{19}$, and $\P^4=\P(\dual{V_5})$,
on
which $\HHH_5$ acts, with eigenspaces of $\iota$ yielding $\P^1_-,\P^2_+
\subseteq\P^4$. We will write the standard generators $\sigma,\tau$
of $\HHH_5$ and $\HHH_{20}$ as $\sigma_5,\tau_5$ and $\sigma_{20},\tau_{20}$
respectively, for clarity. We use coordinates $x_0,\ldots,x_4$
on $\P^4$ and $y_0,\ldots,y_{19}$ on $\P^{19}$, with $\sigma_5(x_i)=
x_{i-1}$, $\tau(x_i)=\xi_5^{-i}x_i$, for $\xi_5$ a fixed primitive
fifth root of unity, and similarly $\sigma_{20}(y_i)=y_{i-1}$
and $\tau_{20}(y_i)=\xi_{20}^{-i}y_i$ for $\xi_{20}$ a fixed 20th
root of unity.

Let $\Delta\subseteq \P^2_+$ be the union of the six lines
$$
\hbox{$x_0=0$ and $x_0+\xi^i x_1+\xi^{-i}x_2=0$, $i\in \boldz_5$,}
$$
for $\xi$ a primitive fifth root of unity. Here we use coordinates
$(x_0,x_1,x_2,x_2,x_1)$ on $\P^2_+$.

In [ADHPR2], \S 8, given a point $a=(a_0,a_1,a_2,a_2,a_1)\in\P^2_+$,
a stable rank 2 reflexive sheaf $\E_a$ is constructed. The precise
details of this construction are not important; rather, only the properties
are, which we summarize here.

\definition{Elines} For a general point $a=(a_0:a_1:a_2)\in\P^2_+$, let 
$$E_{00}:=\lbrace x_1-x_4= x_2-x_3=a_0x_0+2a_1x_1+2a_2x_2=0\rbrace\subset
\P^2_+$$  
and define $E_{ij}:=\sigma_5^i\tau_5^j E_{00}$, for $i,j\in\boldz_5$ 
to be the translates of $E_{00}$ by the Heisenberg group $\HHH_5$.

\proposition{reflexive}
\item{\rm (1)}
$\E_a$ is a stable rank $2$ reflexive sheaf on $\P^4$ 
with Chern classes $c_1=-1$, $c_2=9$, $c_3=25$ and
$c_4=50$. 
\item{\rm (2)}
The singular locus ${\rm Sing}(\E_a)$ of the sheaf $\E_a$ (i.e.,
the locus where it is not locally free)
consists  of the $25$ disjoint lines $E_{ij}$, $i,j\in\boldz_5$.
Moreover, $\E_a(3)$ is generated by sections outside the union of
the lines $E_{ij}$.
\item{\rm (3)} $h^0(\E_a(3))=2$, and the zero locus of a general
section $s\in H^0(\E_a(3))$ is a smooth non-minimal abelian surface
$\tilde A_{a,s}\subset\P^4$, of degree 15, sectional genus 21,
invariant under $\HHH^e_5=\HHH_5\rtimes\langle \iota\rangle$. 
The $25$ lines $E_{ij}$
are exceptional lines on $\tilde A_{a,s}$, 
and they form the canonical divisor of $\tilde A_{a,s}$.
Finally, the polarization $\L$ on the minimal model of $\tilde A_{a,s}$ 
is of type $(1,20)$.

\proof: See [ADHPR2], \S 8. 
\Box

\lemma{H20lemma} 
For general $a\in \P^2_+$, $s\in H^0(\E_a(3))$, the action of
$\HHH^e_5$ on $\P^4$ yields an action of $\HHH^e_5$ on $\tilde A_{a,s}$,
which in turn induces an action on the minimal model $A_{a,s}$ of
$\tilde A_{a,s}$. After a suitable choice of origin in $A_{a,s}$,
the action of $\iota$ on $A_{a,s}$ is given by negation, and the action
of $\HHH_5$ on $A_{a,s}$ is given by translation by the subgroup 
$\boldz_5\times \boldz_5\cong 4K(\L)\subseteq K(\L)$. 

\proof: The action on $\tilde A_{a,s}$ clearly descends to an action
on $A_{a,s}$, and it was already shown in the proof of [ADHPR2], Proposition
8.4, that after a suitable choice of origin, $\iota$ acts as negation.
We continue to use this choice of origin.
Consider next some non-trivial element $\psi\in \HHH_5$, so viewing
$\psi$ as an automorphism of $\P^4$, it has order 5. Furthermore,
it has precisely 5 fixed points in $\P^4$, which form an orbit
under $\HHH_5$. So $\psi$ is clearly
not the identity on $\tilde A_{a,s}$.
Now $\psi$ descends to an automorphism
$\psi:A_{a,s}\rightarrow A_{a,s}$ of order 5. We first need to show
that $\psi$ has no fixed points. However, clearly none of the lines
$E_{ij}$ contains a fixed point of $\psi$, and by \ref{reflexive}, (2),
we can choose $s$ suffiently general so that $A_{a,s}$ does not
contain any given finite set of points disjoint from $\bigcup E_{ij}$. 
Thus the action of $\psi$ on $A_{a,s}$ for $a,s$ general is fixed-point
free.

Now by [BL], (1.1),  
there is an $x\in A_{a,s}$ such that if $t_x$ denotes translation
by $x$, then $\psi'=t_x\circ \psi$ is a group automorphism of $A_{a,s}$,
and the number of fixed points of $\psi'$ and the number of fixed
points of $\psi$ coincide (provided we count an infinite number of fixed
points as being zero). Thus, since $\psi$ has no fixed points, we know
that $\psi'$ must have an infinite number of fixed points. There
are two cases. Either $\psi'$ is the identity, or $\psi'$ is the identity
on an elliptic curve $E\subseteq A_{a,s}$. This gives an exact sequence
$$
\exact{E}{A_{a,s}}{E'}
$$
and $\psi'$ descends to an automorphism on the elliptic curve $E'$.
However, there is no elliptic curve with an automorphism of order $5$,
so this case is ruled out. Thus $\psi'$ is the identity, and $\psi$
is a translation.

Since the only translations which preserve the line bundle $\L$ are 
translations by elements of $K(\L)$, by definition of $K(\L)$, we obtain the
desired result.
\Box 

\theorem{main120theorem} 
\item{\rm (1)} Let $A\subseteq \P^{19}$ be a general $\HHH_{20}$-invariant
abelian surface embedded via a symmetric line bundle $\L$. Then for any
odd 2-torsion point $z\in A\cap \P^-_8$, the subspace 
$V\subseteq H^0(A,\L)$ of sections vanishing on an $\HHH_5$-orbit of $z$
(using the canonical inclusion $\HHH_5\subseteq\HHH_{20}$ given
by $\sigma_5\mapsto \sigma_{20}^4$, $\tau_5\mapsto\tau_{20}^4$) 
is five-dimensional,
and is isomorphic to $V_5$ as an $\HHH_5$-representation. This yields
a canonical identification of
$\P^4=\P(\dual{V_5})$ and $\P(\dual{V})$. Furthermore, the rational map
$A\rto \P(\dual{V})$ induced by the subspace $V$ lifts to an embedding
$\tilde A\rightarrow \P(\dual{V})$, where $\tilde A$ is the blow-up
of $A$ along the $\HHH_5$-orbit of $z$. Finally, the image of $\tilde A$
is the zero-set of some section $s\in H^0(\P(\dual{V}),\E_a(3))$, for
some $a\in\P^2_+$.
\item{\rm (2)} If $A,A'\subseteq\P^{19}$ are two such surfaces as in (1),
then the images of $\tilde A$ and $\tilde A'$ in $\P^4\cong
\P(\dual{V})\cong \P(\dual{(V')})$ coincide if and only if
the canonical level structures on $A$ and $A'$ are $H$-equivalent, with
$H=4K(1,20)$.

\proof: Consider a non-singular surface $\tilde A_{a,s}\subseteq
\P^4=
\P(V^{\vee}_5)$ given by \ref{reflexive}. Let $A_{a,s}$ be the minimal model
of $\tilde A_{a,s}$. The surface $A_{a,s}$ 
carries a degree $(1,20)$ line bundle $\L$ with kernel $K(\L)$. 
We can assume $\L$ is symmetric. Note that by \ref{H20lemma}, the
action of $\HHH_5$ on $\P(V^{\vee}_5)$ induces a symplectic
identification of the group $4K(\L)$ with $K(1,5)$. 
Let $W_{20}= H^0(A_{a,s},\L)$. After choosing a symplectic identification
of $K(\L)$ with $K(1,20)$ compatible with the identification of $4K(\L)$
with $K(1,5)$, we obtain a representation of $\HHH_{20}$ on $W_{20}$
so that the action of $\sigma_5$ on $A_{a,s}$ coming from the embedding
of $\tilde A_{a,s}$ in $\P(V^{\vee}_5)$ coincides with $\sigma_{20}^4$,
and similarly the action of $\tau_5$ coincides with $\tau_{20}^4$. 
Note there is only one such identification up to $H$-equivalence, where
$H=4K(1,20)$, so we obtain a well-defined $H$-level structure on $A_{a,s}$.

Now $\tilde A_{a,s}$ is obtained by choosing a five-dimensional
subspace $W_5\subseteq W_{20}$, and using this subspace to map $A_{a,s}$
birationally into $\P(\dual{W}_5)$, a projective four-space. It follows
however from \ref{H20lemma} that in fact $W_5$ must be an 
$\HHH^e_5$-subrepresentation of $W_{20}$ coming from the natural inclusion
of $\HHH^e_5$ in $\HHH^e_{20}$ given by $\sigma_5\mapsto \sigma_{20}^4$
and $\tau_5\mapsto \tau^4_{20}$. We also have a canonical isomorphism,
up to a scalar multiple, of $V_5$ with $W_5$ as $\HHH^e_5$-representations,
as these are both irreducible representations. Furthermore, the argument
in [ADHPR2], Proposition 8.4 says that in fact $E_{00}\subseteq
\P^2_+$ is the image of the blow-up of an odd two-torsion
point of $A_{a,s}$, and hence all the elements of $W_5$, as
sections of $\L$, must vanish at one of the four odd two-torsion points
of $A_{a,s}$, say $z$. Since $W_5$ is invariant under $\HHH_5$, every 
element of $W_5$ vanishes along the entire $4K(\L)$-orbit of $z$,
which consists of 25 points. Consider the subspace
$W\subseteq W_{20}$ of all sections vanishing on this orbit, so
$W_5\subset W$. Note 
$H^0(\tilde A_{a,s},\O_{\tilde A_{a,s}}(1))\cong W$, and so $\dim W>5$
contradicts linear normality of $\tilde A_{a,s}$. Thus $W=W_5$. 

From this, we see that
the embedding $\tilde A_{a,s}\subseteq \P(V^{\vee}_5)$ is uniquely
determined by the choice of $H$-level structure on $A_{a,s}$ and 
a choice of an odd two-torsion point. However,
given two different choices $z$ and $z'$ of odd two-torsion points
inducing birational maps $\phi,\phi':A_{a,s}\rightarrow \P(\dual{V}_5)$,
there is an $\alpha\in \{\sigma_{20}^{10},\tau_{20}^{10},\sigma_{20}^{10}
\tau_{20}^{10}\}$ with $\alpha(z)=z'$, and then $\phi\circ \alpha=\phi'$.
Thus the surface $\tilde A_{a,s}$ is determined precisely by an $H$-level
structure on $A_{a,s}$. This represents only a finite number of
possible choices. Since the family of surfaces $\tilde A_{a,s}$
parameterized by $a$ and $s$ is three-dimensional, as is the
moduli space of $(1,20)$-polarized abelian surfaces, we conclude
that for general choice of $a,s$, $A_{a,s}$ is general in
the moduli space $\A_{20}$. This, together with the above description
of how $\tilde A_{a,s}$ is obtained from $A_{a,s}$, gives (1) and (2).
\Box

\corollary{A20Hcor} The moduli space  $\A^H_{20}$ is rational.

\proof: \ref{main120theorem} shows that $\A^H_{20}$ is birationally
equivalent to the moduli space of pairs $(a,s)$ where $a\in \PP^2_+$
and $s\in \P(H^0(\P^4,\E_a(3)))$. Now the construction of the sheaf
$\E_a$ given in [ADHPR2] can be done in a relative fashion Zariski locally, 
in the
sense that there is a Zariski open set $U\subseteq \P^2_+$ and
a sheaf $\E$ on $U\times \P^4$ such that
$\E|_{\{a\}\times \P^4}\cong \E_a$ for $a\in U$.
If $p_1,p_2$ are the projections of $\P^2_+\times \P^4$
onto the first and second
factors, then it is clear that $p_{1*}(\E\otimes p_2^*\O_{\PP^4}(3))|_{U'}$ 
is a rank two vector bundle
for some open subset $U'\subseteq U$. Thus $\A^H_{20}$ is birational
to the projectivization of this rank two vector bundle. However, such  
a $\PP^1$-bundle is a rational variety.
\Box

\proposition{quintic}
\item{\rm (1)} Let  $\lbrace s_1, s_2\rbrace$ be a basis of $H^0(\E_a(3))$. 
Then $V_{20,a}=\lbrace s_1\wedge s_2=0\rbrace\subset\P^4$ is 
a quintic hypersurface
containing the pencil of abelian surfaces 
$\{\tilde A_{a,\lambda s_1+\mu s_2}\,|\,(\lambda,\mu)\in \P^1\}$.
The base locus of this pencil is the union of the lines $E_{ij}$.
\item{\rm (2)} For a general choice of the 
parameter $a=(a_0:a_1:a_2)\in\P^2_+$,
$V_{20,a}$ is the unique quintic hypersurface containing the configuration 
of lines $\displaystyle\bigcup_{i,j\in\boldz_5} E_{ij}$. Furthermore, 
the quintic $V_{20,a}$ is $\HHH_5$ invariant.
 
\proof: (1) The fact that $V_{20,a}$ is a quintic hypersurface follows from
the fact that $c_1(\E_a(3))=5$, and this quintic clearly contains the
given pencil of abelian surfaces. The fact that the union of lines $E_{ij}$
is the base locus of the pencil follows from \ref{reflexive}, (2).

(2) We use the idea in [Au], 2.2. The group $\HHH_5$
acts on $H^0(\PP^4,\O_{\PP^4}(5))$ as a sum of characters as follows.
Denote by $V_{r,s}$ the one-dimensional
representation of $\HHH_5$ in which $\sigma$ acts by multiplication by
$\xi^r$ and $\tau$ acts by multiplication by $\xi^s$, where $\xi$ is a 
primitive fifth root of unity. Then
$$H^0(\PP^4,\O_{\PP^4}(5))=6V_{0,0}\oplus \bigoplus_{(r,s)\not=(0,0)}
5V_{r,s}.$$
Indeed, a basis of $\HHH_5$-invariant quintics is given by
{\settabs \+\qquad\qquad$\xi_0=x_0x_1x_2x_3x_4$,&
\qquad\qquad$\xi_1=\sum_{i\in\boldz_5}x_ix_{i+2}^2x_{i+3}^2$,&
\qquad\qquad$\xi_2=\sum_{i\in\boldz_5}x_i^3x_{i+2}x_{i+3}$,\cr
\+$\gamma_0=x_0x_1x_2x_3x_4$,&$\quad\gamma_1=\sum_{i\in\boldz_5}x_ix_{i+2}^2x_{i+3}^2$,&
$\gamma_2=\sum_{i\in\boldz_5}x_i^3x_{i+2}x_{i+3}$,\cr
\medskip
\+$\gamma_3=\sum_{i\in\boldz_5}x_i^3x_{i+1}x_{i+4}$,&\quad$\gamma_4=
\sum_{i\in\boldz_5}x_ix_{i+1}^2x_{i+4}^2$,&
$\gamma_5=\sum_{i\in\boldz_5}x_i^5 - 5x_0x_1x_2x_3x_4$.\cr}
\noindent
A basis for $5V_{r,s}$ is given by
$$
B_{r,s}=\{\sum_{i=0}^4 \xi^{ri}\prod_{j=0}^4 x_{i+j}^{m_j}\,|\,
\sum_{j=0}^4 j m_j\equiv -s\bmod 5, \sum_{j=0}^4 m_j=5\}.
$$
On the other hand, by construction, $H^0(\O_{\cup E_{ij}}(5))$ is
six times the regular representation of $\boldz_5\times\boldz_5$, 
so, by Schur's Lemma,
the restriction map $\rho:H^0(\O_{\P^4}(5))\rightarrow H^0(\O_{\cup E_{ij}}(5))$
decomposes as $\rho=\bigoplus \rho_{r,s}$, where $\rho_{0,0}:6 V_{0,0}
\rightarrow 6V_{0,0}$ and $\rho_{r,s}:5V_{r,s}\rightarrow 6V_{r,s}$ for
$(r,s)\not=(0,0)$. As a consequence, $H^0(\I_{\cup E_{ij}}(5))
=\ker\rho=\bigoplus_{r,s} \ker \rho_{r,s}$.

Thus, in order to prove (2), we check that the mappings $\rho_{r,s}$
are injective for $(r,s)\not=(0,0)$, while $\ker \rho_{0,0}$
is one dimensional. For $(r,s)\not=(0,0)$, this is an open condition on 
$\P^2_+$, so
it suffices to check this in a special case. We take $a_0=1$ and $a_1=a_2=0$.
Using the above bases, this is a straightforward calculation, and we omit
the details.

Now for general $a=(a_0,a_1,a_2)\in \P^2_+$, one can check that the 
$\HHH_5$-invariant quintic
$$\eqalign{ (a_0^5-8a_1^5-8a_2^5+15a_0^3a_1a_2)\gamma_0+(a_0^4a_1+8a_1^3a_2^2-4a_0a_2^4)
\gamma_1+(a_0^3a_2^2-2a_0^2a_1^3-4a_0a_1a_2^3)\gamma_2\cr
+(a_0^3a_1^2-4a_0a_1^3a_2-2a_0^2a_2^3)\gamma_3+
(a_0^4a_2+8a_1^2a_2^3-4a_0a_1^4)\gamma_4+a_0^3a_1a_2\gamma_5=0.}
$$
contains $E_{00}$, and hence $\bigcup_{ij} E_{ij}$ by Heisenberg invariance.
Thus $\ker \rho_{0,0}$ is generically at least dimension one. To show
it is always dimension one is then an open condition, which can again
be checked at $a=(1,0,0)$. Note in this case, $V_{20,a}$ is given by the equation
$x_0x_1x_2x_3x_4=0$. 
\Box

\proposition{100nodes} For a general parameter $a\in\P^2_+$
the quintic $V_{20,a}$ has $100$ ordinary double points as singularities, four of them
on each line $E_{ij}$.

\proof: We only have a computational proof of this fact. By using [BS] or [GS], 
one can easily find choices of $a$ for which $V_{20,a}$ has $100$ ordinary
double points as singularities. We need to show that this is actually
the generic situation.
By computing, again with [BS] or [GS], a standard basis for the Jacobian ideal, one sees
that the singular locus of $V_{20,a}$ 
is supported on $\displaystyle\cup_{i,j\in\boldz_5} E_{ij}$.
Moreover, if $f$ denotes the equation of $V_{20,a}$ then\medskip

\settabs \+\qquad\qquad\qquad${\partial f\over\partial x_2}-{\partial f\over\partial x_3}\in I_{\P^2_+}$&
\qquad\qquad\qquad$a_0\left({\partial f\over\partial x_2}+{\partial f\over\partial x_3}\right)-
2a_2{\partial f\over\partial x_0}\in I_{E_{00}}$\cr
 \+\qquad\qquad\qquad$\displaystyle{\partial f\over\partial x_2}-{\partial f\over\partial x_3}\in I_{\P^2_+}$&
\qquad\qquad\qquad$\displaystyle a_0\left({\partial f\over\partial x_2}+{\partial f\over\partial x_3}\right)-
2a_2{\partial f\over\partial x_0}\in I_{E_{00}}$\cr
\medskip
 \+\qquad\qquad\qquad$\displaystyle{\partial f\over\partial x_1}-{\partial f\over\partial x_4}\in I_{\P^2_+}$&
\qquad\qquad\qquad$\displaystyle a_1\left({\partial f\over\partial x_2}+{\partial f\over\partial x_3}\right)-
2a_2{\partial f\over\partial x_1}\in I_{E_{00}}$\cr
\medskip
\centerline{$\displaystyle a_1{\partial f\over\partial x_0}-a_0{\partial f\over\partial x_1}\in I_{E_{00}}$.}
\noindent
Thus, for $a_0\ne 0$, the singularities of $V_{20,a}$ 
on $E_{00}$ are defined by ${\partial f\over\partial x_0}$.
Restricted to $E_{00}$, 
${\partial f\over\partial x_0}$ has simple roots, and one checks in fact that,
for general choices, these points are $A_1$ singularities on $V_{20,a}$. 
Therefore, by $\HHH_5$ symmetry, 
$V_{20,a}$ has $100$ nodes as singularities. \Box

\proposition{hodge} 
Let $a\in\P^2_+$ be a general parameter and let $\widehat V_{20,a}$ 
be the small resolution of the hypersurface
$V_{20,a}\subset\P^4$ obtained by blowing up one of the smooth abelian surfaces contained in the pencil traced on $V_{20,a}$.
Then
\item{(1)} ${\rm defect\,}(V_{20,a})=1$, while the topological Euler characteristic
$e(\widehat V_{20,a})=200-2(\sharp{\,\rm nodes})=0$
\item{(2)} $\widehat H$, the pullback of a hyperplane section of $V_{20,a}$ and $\widehat A$, the pullback
of an abelian surface in the pencil form a basis for $\Pic(\widehat V_{20,a})/tors
\cong\boldz^2$.
\item{(3)} Finally, $\dim_{\CC} H^1({\widehat V_{20,a}},\Theta_{\widehat V_{20,a}})
=h^1(\Omega_{\widehat V_{20,a}}^2)=2$, hence 
$\P^2_+$ is birational to the moduli of these quintic hypersurfaces.

\proof: The calculation of the defect is easily done on [BS] or [GS] and is
standard; we omit the details. This gives (1). Note the calculation
of the defect implies
$h^{1,2}(\widehat V_{20,a})=2$ from which (3) follows. 
The value of the topological Euler
characteristic gives $h^{1,1}(\widehat V_{20,a})=2$. 

For (2), let $\ell$ be the class in $H^4(\widehat V_{20,a},\boldz)$ 
of the proper transform of $E_{00}$
in $\widehat V_{20,a}$, and let $e$ be the class of one of the exceptional
curves. Then $\widehat H.\ell =1$, $\widehat H. e=0$, 
and $\widehat A.e=-1$. Thus the this intersection matrix is unimodular,
so it shows that $\widehat H,\widehat A$ form a basis of 
$H^2(\widehat V_{20,a},\boldz)/tors\cong \Pic(\widehat V_{20,a})/tors$, 
as desired.
\Box

\remark{finalremarks}
(1) In light of [GPav], it would be interesting to compute the Brauer
group of $\widehat V_{20,a}$. Conjecturally, it would be $\boldz/5\boldz
\times \boldz/5\boldz$.

(2) The Calabi-Yau $\widehat V_{20,a}$ has at least three minimal models;
the first is the one above obtained by blowing up one of the abelian
surfaces. The $100$ exceptional curves can be flopped, giving another
minimal model with a pencil of non-minimal abelian surfaces, with base-locus
the proper transform of $\bigcup_{ij} E_{ij}$. Finally, the curves $E_{ij}$
can be simultaneously flopped, to obtain a minimal model with a base-point-free
pencil of minimal abelian surfaces. We have not, however, completely determined
the K\"ahler cone of $\widehat V_{20,a}$, so there may be some additional
interesting structure here.

\references

\item{[Au]} Aure, A.,
 ``Surfaces on quintic threefolds associated to the Horrocks-Mumford bundle'',
in {\it Lecture Notes in Math.}, {\bf 1399}, (1989), 1--9, Springer.
\item{[ADHPR1]} Aure, A.B., Decker, W., Hulek, K., Popescu, S., Ranestad, K.,
``The Geometry of Bielliptic Surfaces in $\Pfour$'',
{\it Internat. J. Math.}, {\bf 4}, (1993) 873--902.
\item{[ADHPR2]} Aure, A.B., Decker, W., Hulek, K., Popescu, S., Ranestad, K.,
``Syzygies of abelian and bielliptic surfaces in ${\bf P}\sp 4$'',
{\it Internat. J. Math.} {\bf 8}, (1997), no. 7, 849--919.
\item{[BBD]} Bak, A., Bouchard, V., Donagi, R., ``Exploring a new peak
in the heterotic landscape,''
preprint, {\tt arXiv:0811.1242}. 
\item{[BS]} Bayer, D., Stillman, M.,
``Macaulay: A system for computation in
        algebraic geometry and commutative algebra
Source and object code available for Unix and Macintosh
        computers''. Contact the authors, or download from
        {\bf math.harvard.edu} via anonymous ftp.
\item{[Be]} Beauville, A., ``Vari\'et\'es de Prym et jacobiennes
interm\'ediares,'' {\it Ann. Sci. \'Ecole Norm. Sup. (4)}, {\bf 10},
(1977), 309--391.
\item{[BL]} Birkenhake, Ch., Lange, H., ``Fixed-point free automorphisms
of abelian varieties,'' {\it Geom.\ Dedicata},  {\bf 51}, (1994), 201--213.
\item{[BC]} Borisov, L., C\u{a}ld\u{a}raru, A., ``The Pfaffian-Grassmannian
derived equivalence,'' preprint, 2006, {\tt math/0608404}.
\item{[CD]} Candelas, P., Davies, R., ``New Calabi-Yau manifolds with small
Hodge numbers,'' preprint, {\tt arXiv:0809.4681}.
\item{[GS]} Grayson, D., Stillman, M.,
``Macaulay 2: A computer program designed to support
computations in algebraic geometry and computer algebra.''
Source and object code available from
{\tt http://www.math.uiuc.edu/Macaulay2/}.
\item{[GH]} Griffiths, P., Harris, J., {\it Principles of algebraic
geometry,} J. Wiley and sons, New York, 1978. 
\item{[Gri1]} Gritsenko, V., ``Irrationality of the moduli spaces of 
polarized abelian surfaces. With an appendix by the author and K. Hulek''
in Abelian varieties (Egloffstein, 1993), 63--84, de Gruyter, Berlin, 1995.
\item{[Gri2]} Gritsenko, V., ``Irrationality of the moduli spaces of 
polarized abelian surfaces'', {\it Internat. Math. Res. Notices} (1994), 
no. 6, 235 ff., approx. 9 pp. (electronic).  
\item{[GPav]} Gross, M., Pavanelli, S., ``A Calabi-Yau threefold with Brauer 
group $(\boldz/8\boldz)^2$,''
{\it Proc. Amer. Math. Soc.} {\bf 136} (2008), 1--9.
\item{[GP1]} Gross, M., Popescu, S., ``Equations of $(1,d)$-polarized
abelian surfaces,'' {\it Math. Ann.} {\bf 310}, (1998) 333--377.
\item{[GP2]} Gross, M., Popescu, S., ``The moduli space of $(1,11)$-polarized
abelian surfaces is unirational,'', 
{\it Compositio Math.} {\bf 126}, (2001), 1--23.
\item{[GP3]} Gross, M., Popescu, S., ``Calabi-Yau 3-folds and moduli of
abelian surfaces I'', 
{\it Compositio Math.} {\bf 127}, (2001), 169--228.
\item{[HT]} Hori, K., Tong, D., ``Aspects of non-abelian gauge dynamics
in two-dimensional ${\cal N}=(2,2)$ theories,'' {\it J. High Energy
Phys.,} (2007).
\item{[HKW]} Hulek, K., Kahn, C., Weintraub, S., {\it Moduli Spaces of Abelian
Surfaces: Compactification, Degenerations, and Theta Functions},
Walter de Gruyter 1993.
\item{[HS]} Hulek, K., and Sankaran, G., ``The Kodaira dimension of certain
moduli spaces of abelian surfaces,'' {\it Compositio Math.}, {\bf 90},
(1994) 1--35.
\item{[K]} Kuznetsov, A., ``Homological projective duality for Grassmannians
of lines,'' preprint, {\tt arXiv:math/0610957}. 
\item{[LB]} Lange, H., Birkenhake, Ch., {\it Complex abelian varieties},
Springer-Verlag 1992.
\item{[MR]} Melliez, F., Ranestad, K., ``Degenerations of $(1,7)$-polarized
abelian surfaces,'' {\it Math. Scand.}, {\bf 97} (2005), 161--187.
\item{[Mu1]} Mumford, D., ``On the equations defining abelian varieties'',
{\it Inv. Math.}, {\bf 1}, (1966) 287--354.
\item{[Mu2]} Mumford, D., {\it Abelian varieties}, Oxford 
University Press 1974.
\item{[R\o d]} R\o dland, E. A., ``The Pfaffian Calabi-Yau, its mirror,
and their link to the Grassmannian $G(2,7)$,'' {\it Comp. Math.}, {\bf 122},
(2000), 135--149.
\item{[Pop]} Popescu, S., {\it On smooth surfaces of degree $\ge 11$
in the projective fourspace,} Ph.D. Thesis, Saarb\"ucken, 1993.
\item{[SR]} Semple, G., Roth, L., {\it Algebraic Geometry}, Chelsea 1937.

\end